\documentclass[a4paper,9pt]{extarticle}

\usepackage[english]{babel}
\usepackage{ae}
\usepackage{graphicx}
\usepackage{tikz}
\usepackage{amsmath, amsthm, amsfonts, amssymb}
	\numberwithin{equation}{section}
	\theoremstyle{definition}
	\newtheorem{exmp}{Example}
\usepackage{mathtools}
\usepackage{nicefrac}
\usepackage{paralist}
\usepackage[footnotesize, justification=centering, skip=5pt]{caption}
\usepackage{subfig}
\usepackage{setspace}

\providecommand{\abs}[1]{\lvert#1\rvert}
\providecommand{\Abs}[1]{\left\lvert#1\right\rvert}
\providecommand{\norm}[1]{\lVert#1\rVert}
\providecommand{\Norm}[1]{\left\lVert#1\right\rVert}

\DeclareMathOperator*{\TV}{TV}

\newcommand\xqed[1]{
  \leavevmode\unskip\penalty9999 \hbox{}\nobreak\hfill
  \quad\hbox{#1}}
\newcommand\demo{\xqed{$\triangle$}}

\begin{document}

\title{\huge Total Variation Denoising on Hexagonal Grids}
\author{Clemens Kirisits \\
		\small
		Computational Science Center, University of Vienna,\\
		\small
		Nordbergstra{\ss}e 15, 1090 Vienna, Austria}

\date{\today}

\maketitle

\begin{abstract}
This work combines three paradigms of image processing: \begin{inparaenum}[i)]\item the total variation approach to denoising, \item the superior structure of hexagonal lattices, and \item fast and exact graph cut optimization techniques. \end{inparaenum} Although isotropic in theory, numerical implementations of the $BV$ seminorm invariably show anisotropic behaviour. Discretization of the image domain into a hexagonal grid seems perfectly suitable to mitigate this undesirable effect. To this end, we recast the continuous problem as a finite-dimensional one on an arbitrary lattice, before focussing on the comparison of Cartesian and hexagonal structures. Minimization is performed with well-established graph cut algorithms, which are easily adapted to new spatial discretizations. Apart from producing minimizers that are closer in the $\ell^1$ sense to the clean image for sufficiently high degrees of regularization, our experiments suggest that the hexagonal lattice also allows for a more effective reduction of two major drawbacks of existing techniques: metrication artefacts and staircasing. For the sake of practical relevance we address the difficulties that naturally arise when dealing with non-standard images.
\end{abstract}

\section{Introduction}
\label{sec:intro}
\subsection{Total Variation Denoising} \label{sec:tvd}
Denoising is a prominent image processing task, and total variation (TV) regularization is one particular approach. It consists in recovering a noise-free image from perturbed data $f\in L^\alpha(\Omega)$ by finding the solution of an optimization problem, which can be written in general form as
\begin{equation} \label{eq:tvlalpha} \tag{TV$L^\alpha$}
		\min_{u\in BV(\Omega)} \lambda \Norm{u-f}_{L^\alpha}^\alpha + \TV(u).
\end{equation}
We restrict our attention to bounded domains $\Omega \subseteq \mathbb{R}^2$ and only consider scalar functions $f,u\geq0$, which have an immediate interpretation as common greyscale images. Possible minimizers of \eqref{eq:tvlalpha} are required to be of bounded variation, in other words, the space $BV(\Omega)=\{u\in L^1(\Omega): \TV(u) < \infty\}$ is adopted for the modelling of images. For Sobolev functions $u\in W^{1,1}(\Omega)\subset BV(\Omega)$ the total variation equals $\int_\Omega \left|\nabla u\right|$. However, it is usually defined via duality as 
\begin{equation}\label{eq:tv}
		\TV(u) = \sup \left\{\int_\Omega u\,\mathrm{div}g\,\mathrm{d}x: g\in C^\infty_0(\Omega,\mathbb{R}^2),\, \left|g(x)\right|\leq 1,\, x\in\Omega\right\}
\end{equation}
to include discontinuous functions, which is actually one of the main virtues of this denoising technique. Equipped with the norm $\left\|\cdot\right\|_{BV} = \left\|\cdot\right\|_{L^1(\Omega)}+ \TV(\cdot)$ the space of functions of bounded variation is Banach. We refer to \cite{AmbFusPal00,Giu84} for theoretical background on $BV(\Omega)$, to \cite{ChaCasCreNovPoc10,CasChaNov11,SchGraGroHalLen09} for more details on its role in image analysis and to \cite{Kee03} for a concise comparison to other image filters.

Returning to \eqref{eq:tvlalpha}, we only consider the two popular cases $\alpha \in\{1,2\}$, but other choices are possible. For $\alpha = 2$ the minimization problem becomes the unconstrained Rudin-Osher-Fatemi model \cite{RudOshFat92,AcaVog94,ChaLio97}, which was originally designed to remove Gaussian noise. An $L^1$ data term is also widely used and seems to be more suitable for impulsive noise \cite{All92,Nik02,Nik04b,ChaEse05}. From an analytical point of view, one major difference between the two models is that the former is strictly convex, thus having a unique minimizer in $BV(\Omega)$, while the latter is only convex, allowing for multiple minimizers. The weighting parameter $\lambda > 0$ controls the trade-off between data fidelity and degree of regularization. Due to its ability to preserve sharp edges, the TV approach to image analysis is considered a substantial improvement over previous least squares methods using quadratic regularization terms. However, it still suffers from several shortcomings, such as staircasing and loss of texture \cite{Mey01}.

Over the last two decades a great number of different numerical methods for solving \eqref{eq:tvlalpha} have been developed. In the words of \cite{ChaChaYip11}, they can be divided into two classes: \begin{inparaenum}[(i)] \item those that ``optimize'' before discretizing and \item those that do it in reverse order\end{inparaenum}. A common approach that belongs to the first group consists in solving an artificial time evolution of the associated Euler-Lagrange equation with some difference scheme. The second class deals with a discrete version of the original functional, thus solves a finite dimensional minimization problem. Graph cuts, which we chose to employ in this work, belong to the latter.
\subsection{Graph Cuts}
Essentially, graph cut minimization techniques make use of \emph{two} equivalences, both of which were identified several decades ago: first, the one between certain discrete functions and cut functions of $s$-$t$-graphs \cite{Iva65,PicRat75} and, second, that between minimum cuts and maximum flows on such graphs \cite{ForFul56}. In order to minimize the objective function, the edge weights in the graph are chosen so that there is a one-to-one correspondence between costs of cuts and values of the function. Since maximum flows can be found in low-order polynomial time with small constants, solutions can be computed very efficiently \cite{AhuMagOrl93,BoyKol04}. Digital images naturally fit into the graph cuts framework, which is why the computer vision community have exploited their benefits for a wide variety of problems. For a brief introduction to graph cuts see \cite{BoyVek06}. For an in-depth treatment with numerous extensions, applications and a comprehensive bibliography we refer to the recent handbook \cite{BlaKohRot11}.

Concerning the problem of TV denoising, there are several works that employed graph cuts to compute exact minima of discrete versions of \eqref{eq:tvlalpha}
\begin{equation}\label{eq:mrf}
		\lambda \sum_{p}{\abs{u_p-f_p}^\alpha} + \sum_{p, q}{\omega_{pq}\abs{u_p - u_q}}
\end{equation}
in similar ways---see e.g.\ \cite{Hoc01,Cha05,DarSig06a,ChaDar09,GolYin09} (in chronological order). Variables $p$ and $q$ denote sites of a discretization of domain $\Omega$ and $\omega_{pq}$ are weights according to some discrete TV (cf.\ Sec.~\ref{sec:discretetv}). Using a level set formulation, energy \eqref{eq:mrf} is decomposed into a sum of binary functions, which in turn are minimized separately. The fact that, on the one hand, TV admits a complete decoupling of levels of $u$ and that, on the other hand, flow information from one level can be reused for another, makes solving \eqref{eq:tvlalpha} highly efficient. Minimization of the multi-label energy \eqref{eq:mrf} in one large sweep with the method from \cite{Ish03} was reported to be considerably less efficient, in terms of both running time and memory consumption \cite{Cha05}.

Although the TV, as defined in \eqref{eq:tv}, is isotropic, that is it assigns the same values to rotated versions of an image, it is well-known that numerical implementations tend to be anisotropic \cite{CasKunPol99,CasChaNov11}. Due to the two-directional nature of the pixel grid (and thus of the flow network), graph cut approaches are no exception. One way of mitigating this undesirable effect is to consider grids with larger neighbourhood structures. Another possibility is to abandon the square lattice altogether, in favour of a hexagonal one.
\subsection{Hexagonal Grids}\label{sec:hexagonalgrids}

That a rectangular lattice might not be the best choice for the sampling and processing of two-dimensional signals has been recognized at least half a century ago \cite{PetMid62,Mer79}. Since then, researchers who investigated the use of hexagonally arranged pixel configurations almost unanimously conclude that they should be preferred over their more common counterparts for a wide variety of applications, such as edge detection, shape extraction or image reconstruction (see e.g.\ the textbook \cite{MidSiv05}, survey articles \cite{Sta01,HeJia05} or, in chronological order, \cite{KamKam92,YabOga02,ConVilBlu05,FaiPet10a}). The virtues of a hexagonal grid---dense packing, better connectivity, higher angular resolution---have also been acknowledged in the field of morphological image analysis \cite{Soi99,Ser84} as well as for novel imaging frameworks \cite{FaiPet10b}. They can be traced back to a few closely related geometric optimality properties \cite{ConSlo93,Hal01}. Inspired by these and their frequent occurrence in nature (insect and human visual systems, honeycombs) hexagonal arrays are heavily used in the area of optical technology (e.g.\ for large telescope mirrors). Finally, let us remark that hexagonal pixel arrangements are already implemented in robotics vision systems,\footnote{http://centeye.com/} largescale media displays\footnote{http://www.smartslab.co.uk/} and electronic paper \cite{HeiZhoKreRajYan09}.

\begin{figure}[!ht]
  \centering
  \includegraphics[height=0.49\textwidth]{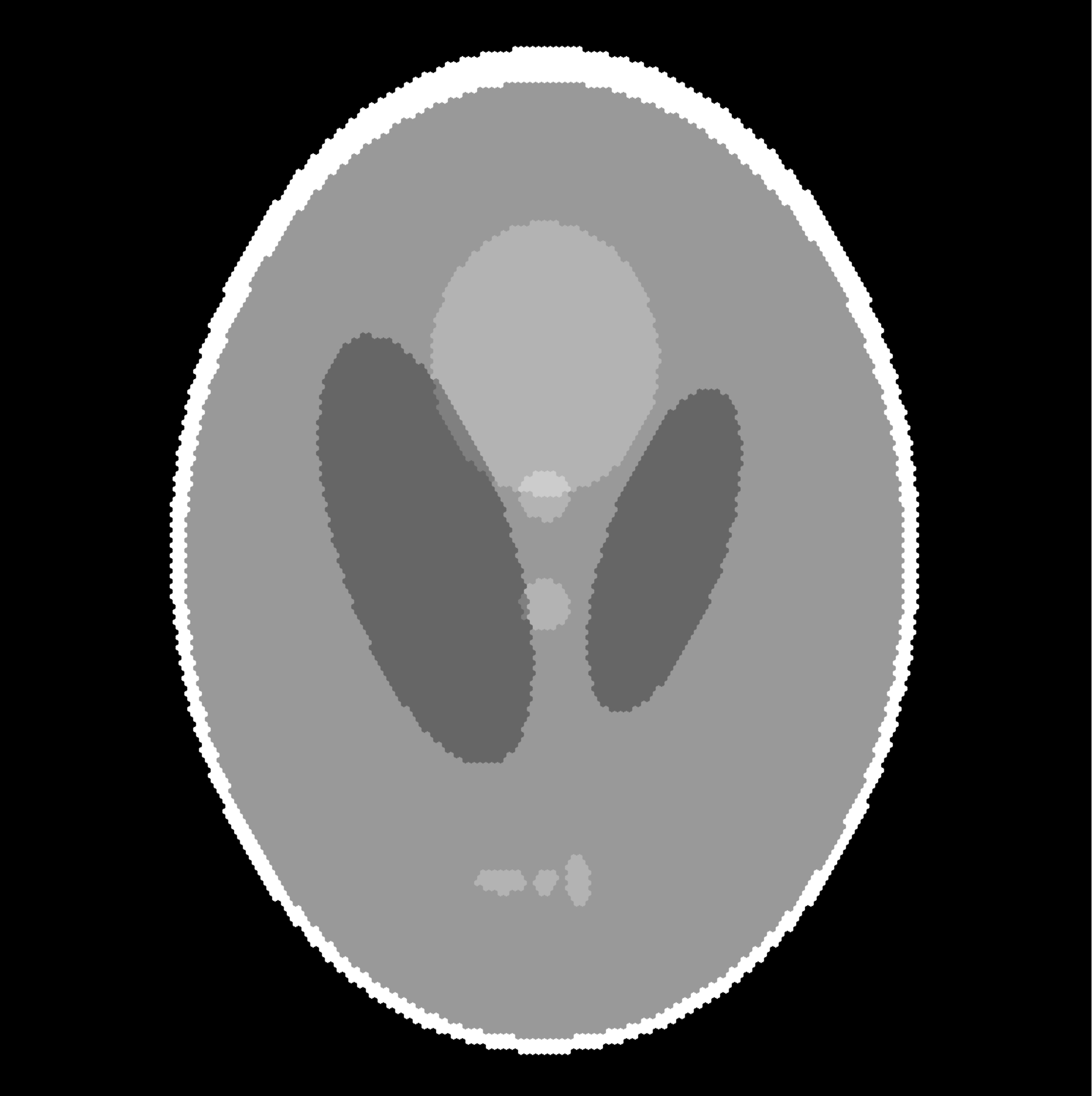}\hfill
  \includegraphics[height=0.49\textwidth]{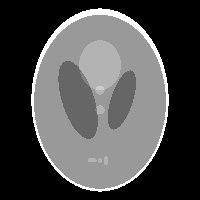}
  \caption{Contrast-enhanced version of the original Shepp-Logan head phantom \cite{SheLog74}, sampled to a $100 \times 100$ square pixel grid (right), and to a hexagonal one with the same resolution (left).}
  \label{fig:shepp-logan}
\end{figure}

\begin{figure}[!ht]
  \centering
  \includegraphics[height=0.49\textwidth]{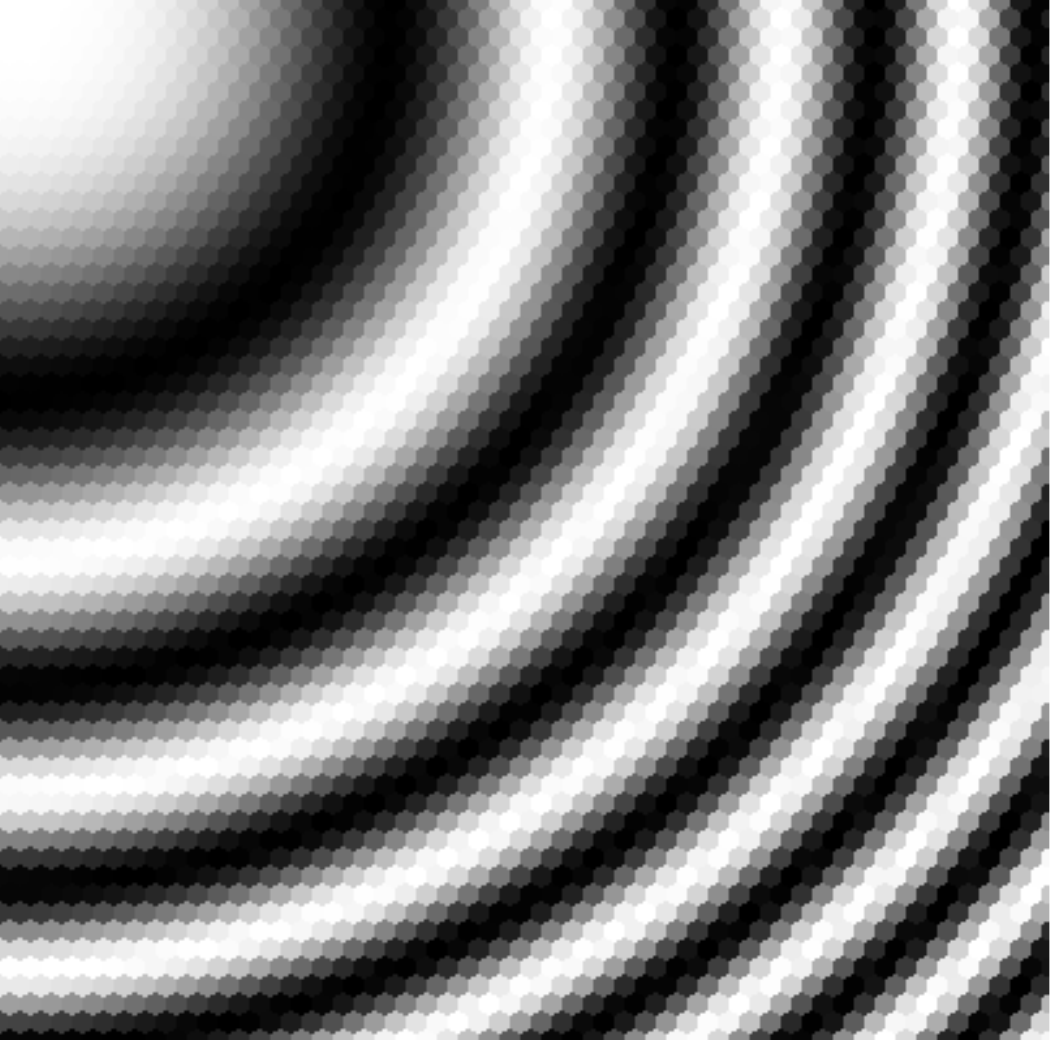}\hfill
  \includegraphics[height=0.49\textwidth]{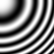}
  \caption{Radially symmetric cosine, $\cos \left(\frac{x^2+y^2}{450}\right)\!$, $(x,y)\in[0,100]\times[-100,0]$, scaled to the interval $[0,255]$ and sampled to both a $54\times 54$ square pixel grid (right) and a hexagonal grid of the same sampling density (left).}
  \label{fig:cos}
\end{figure}

The obvious drawback of a hexagonal grid is its non-alignment with Cartesian coordinates, which means that, in contrast to a square lattice, there is \begin{inparaenum}[a)] no straightforward \item calculus, \item extension to higher dimensions and \item implementability\end{inparaenum}. This, at least to some extent, explains the present lack of hardware devices for efficiently capturing, storing and visualizing hexagonally sampled images, which must be overcome by means of software (see \cite{MidSiv05} for a recently proposed hexagonal image processing framework). In Figs.~\ref{fig:shepp-logan} and \ref{fig:cos} hexagonally sampled images are opposed to common square pixel images. Due to the sufficiently high resolution, the difference between the two sampling schemes seems only marginal in Fig.~\ref{fig:shepp-logan}. In Fig.~\ref{fig:cos}, which features smooth intensity changes, the better visual effects of a hexagonal pixel configuration becomes more apparent. Throughout this note, hexagonal pixels are visualized using so-called \emph{hyperpixels}, i.e.\ groupings of square pixels, that are approximately hexagonal \cite{MidSiv05}.
\subsection{Outline}
When continuous problems such as \eqref{eq:tvlalpha} are brought to finite dimensions, the domain $\Omega$ is usually transformed into a discrete set of points arranged in an orthogonal way. Consequently, functions defined on this lattice inevitably inherit its strongly anisotropic character. However, apart from an easy calculus, there is actually no good reason for choosing a rectangular grid, in particular when isotropy is an important property of the underlying problem as is the case for image processing. The fact that a hexagonal grid features \emph{three} natural directions instead of only two can be expected to mitigate undesirable discretization errors.

To be able to exploit the advantages of the hexagonal grid, we first formulate a finite-dimensional version of \eqref{eq:tvlalpha} on an arbitrary lattice. The corresponding discrete functional is formally derived in Section \ref{sec:derivation}. In Sec.~\ref{sec:lattices} the discretization of domain $\Omega$ is discussed, that is, some background on lattices and neighbourhoods is provided. From that, the definition of discrete images is straightforward (Sec.~\ref{sec:discreteimg}). Sec.~\ref{sec:discretetv} is devoted to the discretization of the TV functional. In Sec.~\ref{sec:discreteproblem} we state our central problem, an ``optimal" finite-dimensional version of problem \eqref{eq:tvlalpha} on an arbitrary lattice, and in Sec.\ \ref{sec:relations} we briefly examine the relations between the continuous and discrete functionals.Finally, in Sec.~\ref{sec:experiments} we present experimental results for hexagonal pixel grid discretizations compared to traditional square ones. We conclude this paper with a brief discussion of our results in Sec.~\ref{sec:conclusion}.

We emphasize that the main contribution of the present work is neither a novel image restoration algorithm nor a theoretical result, but rather a qualitative and purely empirical investigation of different two-dimensional discretization schemes for the problem of TV denoising, focussing on the comparison of square pixel structures with hexagonal ones.
\section{The Discrete Setting} \label{sec:derivation}
\subsection{Planar Lattices}\label{sec:lattices}
A planar lattice is a regular discrete subset of the Euclidean plane. Any basis $B=(b_1,b_2)$ of $\mathbb{R}^2$ generates a lattice $\Lambda(B)$ in the following way
\begin{equation} \label{eq:lattice}
		\Lambda(B) = \left\{ Bp: p \in \mathbb{Z}^2 \right\}.
\end{equation}
Two bases $B_1, B_2$ define the same lattice $\Lambda(B_1) = \Lambda(B_2)$, if and only if one can be obtained from the other with elementary integer column operations, i.e.\ by swapping columns, multiplication of columns by $-1$ and adding integer multiples of one column to another. In other words, there exists a unimodular matrix $U\in\mathbb{R}^{2\times 2}$, an integer matrix with integer inverse, so that $B_1 = B_2U$.
\begin{exmp}
Let us give two examples that will be used throughout this paper. First, if $B=I$ is the identity matrix, $\Lambda(B) = \mathbb{Z}^2$ is the usual unit square lattice. Second, if the basis vectors $b_1,b_2$ have the same length and draw an angle of $\frac{\pi}{3}$ or $\frac{2\pi}{3}$, then $\Lambda(B)$ becomes a regular hexagonal lattice. For the sake of simplicity we only consider horizontally aligned hexagonal lattices, and set accordingly
\begin{equation}\label{eq:hexmatrix}
		H \coloneqq (h_1, h_2) \coloneqq  \begin{pmatrix} \cos 0 & \cos \frac{\pi}{3} \\ \sin 0 & \sin \frac{\pi}{3} \end{pmatrix} = \begin{pmatrix} 1 & \frac{1}{2} \\ 0 & \frac{\sqrt{3}}{2} \end{pmatrix}.
\end{equation}
Fig.~\ref{fig:grid2} shows a horizontally aligned hexagonal lattice.
\begin{figure}
  \centering
  \includegraphics[width=.7\textwidth]{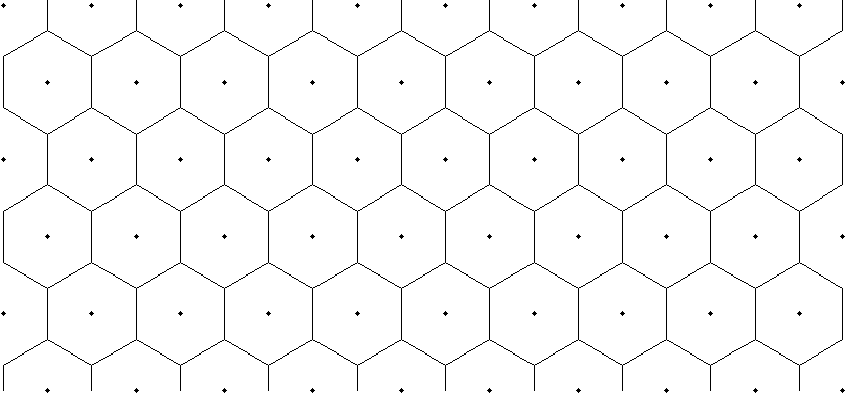}
  \caption{Horizontally aligned regular hexagonal lattice with Voronoi cell tiling.}
  \label{fig:grid2}
\end{figure}
\demo
\end{exmp}
\subsubsection{Lattice Density and Tilings}
From now on, we restrict ourselves to bases $B$ consisting of unit vectors, and introduce an additional parameter $d>0$ controlling the distance between points of lattice $\Lambda(dB)$. A discretization of the image domain $\Omega$ is then defined as
\begin{equation}\label{eq:generalgrid}
		\Omega^d_B \coloneqq \Omega\cap\Lambda(dB).
\end{equation}
We set $N = \abs{\Omega_B^d}$ and, in order to simplify notations, we will identify a lattice point $dBp \in \Omega_B^d$ simply by its ``index" $p$. In general, $d$ does not equal the minimal distance between two lattice points, although it does in the two interesting cases $\Lambda(dI)$ and $\Lambda(dH)$. Parameter $d$ does, however, equal the side length of the lattice's fundamental parallelogram $\Pi(dB) = \left\{ dBx: x\in [0,1)^2 \right\}$. The parallelogram's area $\abs{\det dB}$ is an important quantity, since it is inversely proportional to lattice density: there is one lattice point per parallelogram. From it we can deduce the phenomenon, which is sometimes referred to as the ``denser packing'' of hexagonal lattices:
\begin{equation}\label{eq:density}
		\lvert \det dH \rvert = \frac{\sqrt{3}}{2}d^2 < d^2 = \lvert \det dI \rvert.
\end{equation}
This implies that $\Lambda(dH)$ has a higher density than $\Lambda(dI)$ and that, in general, $\abs{\Omega^d_H} > \abs{\Omega^d_I}$ will hold. Obviously, if we want to compare different lattices, it is important that they are equally dense. Therefore, due to \eqref{eq:density}, we have to shrink the square lattice by a factor of $\sqrt{\nicefrac{\sqrt{3}}{2}}$, or stretch the hexagonal one by $\sqrt{\nicefrac{2}{\sqrt{3}}}$, in order to get a square lattice with the same density as a hexagonal one---and, thus, a discretization of $\Omega$ with approximately the same amount of lattice points.

Associated to every lattice, there is one monohedral tiling, i.e.\ a partition of the plane consisting of congruent shapes, which will be of use later on. It is made up of Voronoi cells (cf.\ Fig.~\ref{fig:grid2}). The Voronoi cell associated to one point $p$ belonging to some discrete set $\Lambda \subset \mathbb{R}^2$ is defined as the set of all points that are closer to $p$ than to any other point from $\Lambda$
\begin{equation}
V(p) = \{ x\in\mathbb{R}^2: \abs{x-p} \leq \abs{x-q},\,\forall q \in \Lambda\setminus\{p\}\}.
\end{equation}
The area of a Voronoi cell associated to lattice $\Lambda(B)$ is equal to that of the fundamental parallelogram $\abs{\det B}$, since Voronoi cells and translations of the fundamental parallelogram both form monohedral tilings of the plane with one lattice point per tile. In the following we will mainly be interested with that part of a Voronoi cell that lies in $\Omega$, i.e.\ with $V_p \coloneqq V(p)\cap\Omega$. For more details on (hexagonal) lattices we refer to \cite{ConSlo93,MidSiv05}.
\subsubsection{Neighbourhoods}
There are several equivalent ways of defining a neighbourhood relation between lattice points. We define it (as in \cite{BoyKol03}) by a set of vectors pointing in the direction of neighbouring sites
\begin{equation} \label{eq:neighborhood}
		\mathcal{N}^\nu_B = \{v_i \in \mathbb{R}^2:i=1, \dotsc, \nu\}.
\end{equation}
These vectors are the ``directions'' of the lattice and, since neighbourship is a symmetric relation, by $v_i$ we actually mean $\pm v_i$. It is also clear that these vectors are linear combinations of the lattice's basis vectors. Two sites $p$ and $q$ are neighbours according to $\mathcal{N}^\nu_B$, if $dB(p-q)$ equals some $v_i \in \mathcal{N}^\nu_B$. We will indicate this with the shorthand $p \sim q$, when there is no confusion about the neighbourhood system.

\begin{figure}[!ht]
	\centering
	\subfloat[Neighbours of the central black point according to $\mathcal{N}^4_I$ (thick lines), $\mathcal{N}^8_I$ (thick and fine) and $\mathcal{N}^{16}_I$ (thick, fine and dashed).]{\label{fig:sneighbors}
	\begin{tikzpicture}
	  \tikzstyle{vertex}=[shape=circle,minimum size=8pt,draw]
	  \tikzstyle{nearest}=[shape=circle,fill=black!60,minimum size=8pt,draw]
	  \tikzstyle{nonnearest}=[shape=circle,fill=black!30,minimum size=8pt,draw]
	  \tikzstyle{farthest}=[shape=circle,fill=black!10,minimum size=8pt,draw]
	  \node (c) at (0,0) [shape=circle, fill=black, minimum size=8pt,draw] {};
	  \node[nearest] (t) at (0,1) {};
	  \node[nearest] (r) at (1,0) {};
	  \node[nearest] (b) at (0,-1) {};
	  \node[nearest] (l) at (-1,0) {};
	  \node[nonnearest] (tr) at (1,1) {};
	  \node[nonnearest] (br) at (1,-1) {};
	  \node[nonnearest] (bl) at (-1,-1) {};
	  \node[nonnearest] (tl) at (-1,1) {};
	  \node[farthest] (ttr) at (1,2) {};
	  \node[farthest] (trr) at (2,1) {};
	  \node[farthest] (rrb) at (2,-1) {};
	  \node[farthest] (rbb) at (1,-2) {};
	  \node[farthest] (lbb) at (-1,-2) {};
	  \node[farthest] (llb) at (-2,-1) {};
	  \node[farthest] (llt) at (-2,1) {};
	  \node[farthest] (ltt) at (-1,2) {};
	  \node[vertex] at (2,0) {};
	  \node[vertex] at (2,2) {};
	  \node[vertex] at (0,2) {};
	  \node[vertex] at (-2,2) {};
	  \node[vertex] at (-2,0) {};
	  \node[vertex] at (-2,-2) {};
	  \node[vertex] at (0,-2) {};
	  \node[vertex] at (2,-2) {};
	  \draw [very thick] (c) -- (t);
	  \draw [very thick] (c) -- (r);
	  \draw [very thick] (c) -- (b);
	  \draw [very thick] (c) -- (l);
	  \draw (c) -- (tr);
	  \draw (c) -- (br);
	  \draw (c) -- (bl);
	  \draw (c) -- (tl);
	  \draw [dashed] (c) -- (ttr);
	  \draw [dashed] (c) -- (trr);
	  \draw [dashed] (c) -- (rrb);
	  \draw [dashed] (c) -- (rbb);
	  \draw [dashed] (c) -- (lbb);
	  \draw [dashed] (c) -- (llb);
	  \draw [dashed] (c) -- (llt);
	  \draw [dashed] (c) -- (ltt);
	\end{tikzpicture}}
	\hspace{30pt}
	\subfloat[Neighbours of the central black point according to $\mathcal{N}^6_H$ (thick lines), $\mathcal{N}^{12}_H$ (thick and fine)]{\label{fig:hneighbors}
	\begin{tikzpicture}
	  \tikzstyle{vertex}=[shape=circle,minimum size=8pt,draw]
	  \tikzstyle{nearest}=[shape=circle,fill=black!50,minimum size=8pt,draw]
	  \tikzstyle{nonnearest}=[shape=circle,fill=black!20,minimum size=8pt,draw]
	  \node (c) at (0,0) [shape=circle, fill=black, minimum size=8pt, draw] {};
	  \node[nearest] (0) at (1*1.155,0) {};
	  \node[nearest] (pi3) at (0.5*1.155,1) {};
	  \node[nearest] (2pi3) at (-0.5*1.155,1) {};
	  \node[nearest] (pi) at (-1*1.155,0) {}; 
	  \node[nearest] (4pi3) at (-0.5*1.155,-1) {};
	  \node[nearest] (5pi3) at (0.5*1.155,-1) {};	  
	  \node[nonnearest] (pi6) at (1.5*1.155,1) {};  
	  \node[nonnearest] (pi2) at (0,2) {};
	  \node[nonnearest] (5pi6) at (-1.5*1.155,1) {};
	  \node[nonnearest] (7pi6) at (-1.5*1.155,-1) {}; 
	  \node[nonnearest] (3pi2) at (0,-2) {};
	  \node[nonnearest] (11pi6) at (1.5*1.155,-1) {};
	  \node[vertex] at (2*1.155,0) {};  
	  \node[vertex] at (1*1.155,2) {};
	  \node[vertex] at (-1*1.155,2) {};
	  \node[vertex] at (-2*1.155,0) {}; 
	  \node[vertex] at (-1*1.155,-2) {};
	  \node[vertex] at (1*1.155,-2) {};
	  \node[vertex] at (2*1.155,2) {};  
	  \node[vertex] at (-2*1.155,2) {};
	  \node[vertex] at (-2*1.155,-2) {};
	  \node[vertex] at (2*1.155,-2) {}; 
	  \draw [very thick] (c) -- (0);
	  \draw [very thick] (c) -- (pi3);
	  \draw [very thick] (c) -- (2pi3);
	  \draw [very thick] (c) -- (pi);
	  \draw [very thick] (c) -- (4pi3);
	  \draw [very thick] (c) -- (5pi3);
	  \draw (c) -- (pi6);
	  \draw (c) -- (pi2);
	  \draw (c) -- (5pi6);
	  \draw (c) -- (7pi6);
	  \draw (c) -- (3pi2);
	  \draw (c) -- (11pi6);
	\end{tikzpicture}}
	\caption{Square lattice (left) and hexagonal lattice (right) with neighbourhoods from \eqref{eq:squareneighbors} and \eqref{eq:hexneighbors}, respectively.}
	\label{fig:neighbors}
\end{figure}
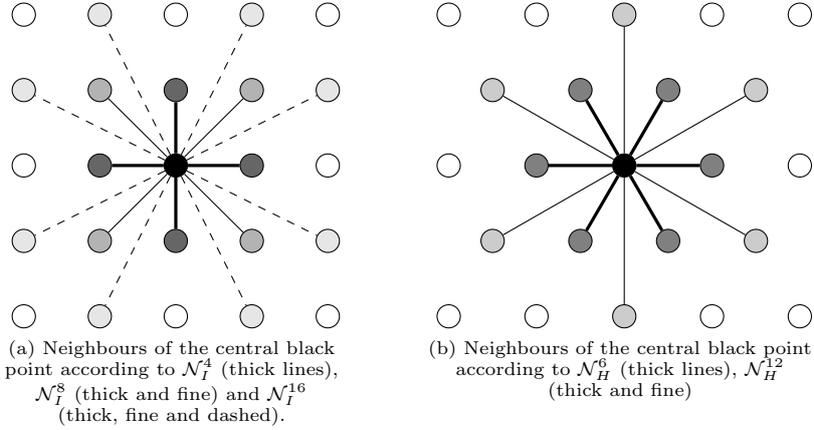

For the square lattice at least three neighbourhood systems (Fig.~\ref{fig:sneighbors}) have been proposed and already used in the context of TV minimization \cite{GolYin09}:
\begin{equation}\label{eq:squareneighbors}
		\begin{aligned}
				\mathcal{N}^4_I		&= \left\{ e_1, e_2\right\},\\
				\mathcal{N}^8_I		&= \mathcal{N}^4_I \cup \left\{ e_1+e_2,e_2-e_1\right\},\\
				\mathcal{N}^{16}_I	&= \mathcal{N}^8_I \cup \left\{ 2e_1+e_2, e_1+2e_2, 2e_2-e_1, e_2-2e_1 \right\}.
  	\end{aligned}
\end{equation}
Here, $e_i$ are the canonical basis vectors and we set $d=1$ for simplicity. For the hexagonal lattice two natural neighbourhoods are (Fig.~\ref{fig:hneighbors}):
\begin{equation}\label{eq:hexneighbors}
		\begin{aligned}
				\mathcal{N}^6_H &= \left\{ h_1, h_2, h_2-h_1\right\}, \\
				\mathcal{N}^{12}_H &= \mathcal{N}^6_H \cup \left\{ h_1+h_2, 2h_2-h_1, h_2-2h_1\right\}.
  	\end{aligned}
\end{equation}
The important point here is that $\abs{h_2-h_1}=\abs{h_1}=\abs{h_2}$. A regular hexagonal lattice is the only planar lattice that permits each of its points to have six nearest neighbours instead of only four or even less.
\subsection{Generalized Discrete Images and their Levels}
\label{sec:discreteimg}
Throughout this note we regard a \emph{generalized discrete image} $\mathbf{u}$ to be a quantized and sampled version of $u \in BV(\Omega)$. That is, $u$ is rounded to a discrete set of grey values $\mathcal{L} \subset \mathbb{R}$, and evaluated at sites of lattice $\Omega_B^d$. According to the notation introduced in the previous section, we write
\begin{equation}\label{eq:gdi}
	\mathbf{u} \coloneqq (u_p)_{p\in\Omega_B^d} \in \mathcal{L}^N, \quad u_p \coloneqq \lfloor u(p) \rfloor,
\end{equation}
where $\lfloor \cdot \rfloor$ stands for quantization.\footnote{Since $u \in BV(\Omega) \subset L^1(\Omega)$, we think of $u(p)$ as $\lim_{r \rightarrow 0^+} (r^2\pi)^{-1} \int_{\mathcal{B}_r(p)}u$, where $\mathcal{B}_r(p)$ is the disk of radius $r$ centred at $p$.} Such an image is ``general" in the sense that for a rectangular domain $\Omega$ and $B=I$ it reduces to the widely used matrix representation of digital images. Additionally, in order to simplify notations, $\mathcal{L}$ is assumed to be a set of $L$ consecutive integers starting from zero, i.e.\ $\mathcal{L}=\{0,\ldots,L-1\}$, where $0$ and $L-1$ are visualized as black and white, respectively, and values in between are represented by different shades of grey (cf.~\ref{fig:cos}).

For $l \in \mathbb{R}$ the super level sets of $u$ are defined as
\begin{equation}\label{eq:levelset}
	E_l = \{x\in\Omega: u(x) > l\}.
\end{equation}
We will denote the characteristic function of a super level set $E_l$ of $u$ simply by $\chi^l$, i.e.\ $\chi^l(x)=1$, if $u(x)>l$ and $\chi^l(x)=0$ otherwise. Image levels are analogously defined in the discrete setting. To each discrete image $\mathbf{u}$, $L$ binary images $\pmb{\chi}^l = (\chi^l_p)_{p\in\Omega_B^d}$, $l \in \mathcal{L},$ are associated, where $\chi^l_p$ equals 1, if $u_p$ is above level $l$, and 0 if $u_p$ is below or equal to $l$. Note that there are actually only $L-1$ relevant levels, since $\chi^{L-1}\equiv 0$ contains no information. This family of binary images is monotone decreasing, i.e.\ it satisfies
\begin{equation} \label{eq:monotony}
		l,k \in \mathcal{L},\quad l > k \quad \Rightarrow \quad \chi^l_p \leq \chi^k_p, \quad p\in\Omega_B^d.
\end{equation}
In other words, for each lattice point $p$ they form a monotone decreasing sequence as depicted in Tab.~\ref{tab:levels}.
\begin{table}
	\centering
	\caption{Super level set characteristics of a discrete image.}
	\label{tab:levels} 
	\begin{tabular}{r *{7}c}
		\hline\noalign{\smallskip}
		$l:$	& 	0	&	1	&	$\cdots$	&	$u_p-1$	&	$u_p$	&	$\cdots$	&	$L-1$	\\
		\noalign{\smallskip}\hline\noalign{\smallskip}
		$\chi^l_p:$	&	1	&	1	&	$\cdots$	&	1		&	0		&	$\cdots$	&	0 \\
		\noalign{\smallskip}\hline
	\end{tabular}
\end{table}
Knowing all levels of an image is equivalent to knowing the image itself. Thus $\mathbf{u}$ can be reconstructed from its levels by $u_p = \sum_l\chi^l_p = \min \{l: \chi^l_p=0\}$.
\subsection{Discrete Total Variation}\label{sec:discretetv}

In this section we demonstrate how to discretize the TV functional in a reasonable way, by combining the coarea and Cauchy-Crofton formulae. While the former states that the TV can be computed by summing level contours, the latter provides a sound way of estimating those contours from function values on a lattice. This is the common approach for graph cut approaches to TV denoising \cite{DarSig06a,GolYin09}. Below, we demonstrate that this approach is easily generalized to arbitrary lattices and, by doing so, make a further point in favour of the hexagonal grid.

For $u\in BV(\Omega)$ the \emph{coarea formula} reads
\begin{equation} \label{eq:coarea}
				\TV(u)	= \int_{-\infty}^\infty \TV(\chi^l)\,\mathrm{d}l,
\end{equation}
where $\chi^l$ is the characteristic function of the set $E_l$ as defined in \eqref{eq:levelset} \cite{AmbFusPal00,Giu84}. The quantity $\mathrm{Per}(E_l)\coloneqq\TV(\chi^l)$ is finite for almost all $l\in\mathbb{R}$ and is commonly called perimeter of the level set. The coarea formula thus states that the TV of a $BV$-function can be computed by summing the perimeters of all its super level sets.\footnote{In the case where $E_l$ has a boundary of class $C^1$ (up to a set which is negligible with respect to the one-dimensional Hausdorff measure) one has $\abs{\partial E_l}=\mathrm{Per}(E_l)$, that is, the perimeter of $E_l$ indeed equals the length of its boundary. In general, however, for $u\in BV(\Omega)$ the topological boundary is not the appropriate object of study and must be replaced with the (measure-theoretic) reduced boundary of $E_l$, a dense subset of $\partial E_l$ \cite{AmbFusPal00,Giu84}.} In the following we denote the $n$-dimensional Hausdorff measure of sets by vertical bars. From the context it should be clear whether the number of elements ($n=0$), length ($n=1$) or area ($n=2$) is meant.
\subsubsection{The Cauchy-Crofton Formula}\label{sec:cutmetric}
The \emph{Cauchy-Crofton formula} is an elegant result from integral geometry, relating the length of a plane curve to the measure of straight lines intersecting it \cite{Car76}. Inspired by the work of \cite{BoyKol03}, we utilize it to estimate from $\mathbf{u}$ the level set perimeters of $u$ in a sound way.

Consider a straight line $k$ in $\mathbb{R}^2$. By identifying it with its \emph{signed} distance from the origin $\rho$ and the angle $\phi$ which it draws with the $x$-axis, line $k$ can be regarded as a point in the set $\mathcal{K} = \mathbb{R}\times [0,\pi)$.\footnote{To make the assignment unique, we set $(\rho,\frac{\pi}{2})$ to be the line $x=\rho$. Note that the Cauchy-Crofton formula is usually stated using \emph{normal parameters} $(p, \varphi)$, i.e.\ its unsigned distance $p$ to the origin and the angle $\varphi$ which is drawn by the $x$-axis and a line normal to $k$ \cite{Car76}.} Let $C$ be a rectifiable plane curve $C$ and $n_C$ the function mapping each line to the number of times it intersects with $C$. The Cauchy-Crofton formula states that
\begin{equation} \label{eq:crofton}
	\abs{C} = \frac{1}{2}\int_\mathcal{K} n_C(\rho,\phi) \, \mathrm{d}\rho \, \mathrm{d}\phi .
\end{equation}

One of the virtues of equation \eqref{eq:crofton} is that it can also be used to approximate the length of $C$: by replacing the double integral on the right-hand side with a partial sum, which corresponds to equipping the plane with countably many lines, a reasonable estimate for the curve length can be obtained. More precisely, we have $m$ families of lines $F_i$, each characterized by its angle $\phi_i$ and inter-line distance $\Delta\rho_i$, i.e.
\begin{equation} \label{eq:family}
F_i(\phi_i,\Delta\rho_i) = \{(j\Delta\rho_i, \phi_i): j\in\mathbb{Z}\},\quad i=1,\ldots,m.
\end{equation}
In addition, we assume $0 = \phi_1 < \phi_2 < \cdots < \phi_m < \pi$ and set $\Delta\phi_i = \phi_{i+1}-\phi_i$ for $i \leq m-1$, and $\Delta\phi_m = \pi-\phi_m$. Then, from \eqref{eq:crofton} we get
\begin{equation} \label{eq:crofton-approx}
	\abs{C} \approx \frac{1}{2} \sum_{i=1}^m \sum_{j\in\mathbb{Z}} n_C(j\Delta\rho_i,\phi_i)\Delta\rho_i\Delta\phi_i = \sum_{i,j} n_C(j\Delta\rho_i,\phi_i) \omega_i.
\end{equation}
Here, we used the shorthand $\omega_i = \Delta\rho_i\Delta\phi_i/2$. This approximation formula tells us for each family $F_i$ to count its intersections with $C$ and to multiply the result with half the area of the rectangle in $\mathcal{K}$ corresponding to $F_i$. Moreover, \eqref{eq:crofton-approx} is a reasonable approximation, since it converges to the true length of $C$ as $m\rightarrow\infty$ and $\sup_i\Delta\rho_i \rightarrow 0$, $\sup_i\Delta\phi_i \rightarrow 0$ \cite{Car76}.

Consider a grid graph embedded in the same plane where $C$ lies, consisting of a set of nodes $\Lambda(dB)$ and edges according to some neighbourhood $\mathcal{N}_B^\nu$. Such a graph naturally induces $m=\frac{\nu}{2}$ families of straight lines $F_i$. Furthermore, we assign weights $\omega_i$ from \eqref{eq:crofton-approx} to edges, depending on the family to which they belong. Now, by summing up the weights of all those edges that have been intersected by $C$, we get exactly the right hand side of \eqref{eq:crofton-approx}, provided the lattice is sufficiently fine so that each edge is intersected by $C$ at most once. Due to its higher angular resolution, a hexagonal grid graph can be expected to yield a better approximation than a square one, since it naturally equips the plane with more families of lines.
\begin{exmp}[Weights for $\mathcal{N}_H^6$, $\mathcal{N}_H^{12}$]
In the case of a regular hexagonal lattice $\Lambda(dH)$ with neighbour distance $d>0$ and neighbourhood $\mathcal{N}_H^6$, we have $m=3$ families of lines at angles $\phi_i = (i-1)\frac{\pi}{3}$ and line distance $\Delta\rho_i = \Delta\rho = d\sqrt{3}/2$. Therefore, we get $\omega_i = \omega = \sqrt{3}\frac{\pi}{12}d$.

For $\mathcal{N}_H^{12}$ (cf.\ Fig.~\ref{fig:hneighbors}) we have $m=6$ families at angles $\phi_i = (i-1)\frac{\pi}{6}$ and line distances
\begin{equation} \notag
		\Delta\rho_i = 
		\begin{cases}
				\frac{\sqrt{3}}{2}d,	& \text{for } i\in\{1,3,5\} \\
				\frac{1}{2}d,			& \text{for } i\in\{2,4,6\}\text{.}
		\end{cases}
\end{equation}
Accordingly, we get two different weights: one for edges between nearest neighbours standing at distance $d$, and one for edges connecting sites at distance $d\sqrt{3}$
\begin{equation} \notag
		\omega_i = 
		\begin{cases}
				\frac{\sqrt{3}\pi}{24}d,	& \text{for } i\in\{1,3,5\} \\
				\frac{\pi}{24}d,			& \text{for } i\in\{2,4,6\}\text{.}
		\end{cases}
\end{equation}
The weights for $\mathcal{N}^4_I$, $\mathcal{N}^8_I$ and $\mathcal{N}^{16}_I$ can be found in \cite{GolYin09}.
\demo
\end{exmp}
\subsubsection{TV Estimation on Arbitrary Lattices}\label{sec:tv-est}
 Consider a continuous signal $u\in BV(\Omega)$, which is known to us only at sampling points $\Omega_B^d$, i.e.\ as a discrete image $\mathbf{u}$ on an arbitrary lattice (cf.\ Sec.~\ref{sec:discreteimg}). Clearly, since we do not know $u$, we do not know its contour lines $\partial E_l$ either. We do know, however, the discrete level characteristics $\pmb{\chi}^l$, which tell us that, if $\abs{\chi^l_p-\chi^l_q}=1$ for neighbouring sites $p$ and $q$, then the contour has to lie somewhere between them. In other words, we regard an edge $\{p,q\}$ of the underlying grid graph as intersected by the unknown contour wherever $\chi^l_p\neq\chi^l_q$. Therefore, formula \eqref{eq:crofton-approx} can be written as
\begin{equation} \label{eq:cutmetric}
		\abs{\partial E_l} \approx \sum_{p \sim q} \omega_{pq} \abs{\chi^l_p - \chi^l_q},
\end{equation}
where $\omega_{pq}$ equals $\omega_i$ from \eqref{eq:crofton-approx}, if the edge $\{p,q\}$ is parallel to family $F_i$. Finally, summing \eqref{eq:cutmetric} over all levels $l \in\mathcal{L}$ and using the coarea formula \eqref{eq:coarea} gives an estimate for the TV of the continuous signal:
\begin{equation} \label{eq:cutmetrictv}
		\TV(u) \approx \sum_{l=0}^{L-2}\sum_{p \sim q} \omega_{pq} \lvert \chi^l_p - \chi^l_q \rvert = \sum_{p \sim q} \omega_{pq} \lvert u_p - u_q \rvert \eqqcolon \TV(\mathbf{u}).
\end{equation}
The last equality in follows from the fact that $\sum_l \abs{\chi^l_p - \chi^l_q} = \abs{\sum_l ( \chi^l_p - \chi^l_q )} = \abs{u_p-u_q}$. Of course, the quality of this approximation depends not only on the goodness of the Cauchy-Crofton approximation \eqref{eq:crofton-approx}, i.e.\ on neighbour distance $d$ and (the size $\nu$ of) the neighbourhood, but also on the discretization $\mathcal{L}$ of $\mathbb{R}$. For a relation between the quantized and unquantized discrete TV$L^2$ problems see \cite{ChaDar09}.
\subsubsection{Related Approaches}\label{sec:related}
The authors of \cite{ChaDar09} define a discrete TV, as \emph{any} function $J:\mathbb{R}^N \rightarrow [0,\infty]$ that satisfies a discrete coarea formula
\begin{equation} \label{eq:discretecoarea}
		J(\mathbf{u}) = \int_{\mathbb{R}} J(\pmb{\chi}^l)\,\mathrm{d}l,
\end{equation}
and, accordingly, a perimeter for discrete sets $A\subset \Omega_B^d$ as $\mathrm{Per}(A) = J(\pmb{\chi}^A)$. This definition has two interesting implications. First, such a $J$ has several desirable properties, like positive homogeneity, lower semicontinuity and submodularity. Second, the straightforward discretization of the TV on a square grid
\begin{equation} \label{eq:squarel2tv}
				\int_{\Omega} \lvert\nabla u\rvert \approx \sum_{i,j} \sqrt{(u_{i+1,j}-u_{i,j})^2 + (u_{i,j+1}-u_{i,j})^2}
\end{equation}
does \emph{not} fall into that category, except for the case $\mathcal{L} = \mathbb{B} =  \{0,1\}$ where equation \eqref{eq:discretecoarea} becomes trivial, since $\pmb{\chi}^0 = \mathbf{u}$. The TV discretization on the right hand side of \eqref{eq:squarel2tv} is actually only submodular for $\mathcal{L} = \mathbb{B}$, thus rendering it unsuitable for standard graph cut approaches.

In any case, our derivation of a discrete TV does fulfil \eqref{eq:discretecoarea}---note that we actually already used the discrete coarea formula in \eqref{eq:cutmetrictv}---, and infinitely many more variants can be constructed, namely by using arbitrary neighbourhoods with arbitrary positive weights. However, only \eqref{eq:cutmetrictv} is justified by the Cauchy-Crofton formula.
\subsection{A General Version of the Discrete Problem}\label{sec:discreteproblem}
After discretizing the fidelity term of \eqref{eq:tvlalpha} in a straightforward way
\begin{equation}\label{eq:discretefidelity}
		\norm{u-f}_{L^\alpha(\Omega)}^\alpha =
		\int_\Omega \abs{u-f}^\alpha\,\mathrm{d}x \approx
		\sum_{p}\abs{V_p}\abs{u_p-f_p}^\alpha,
\end{equation}
we finally arrive at a discrete version of the TV$L^\alpha$ functional.
\begin{equation}\label{eq:discrete-functional}
		F_\alpha(\mathbf{u}) \coloneqq \lambda \sum_{p}\abs{V_p}\abs{u_p-f_p}^\alpha + \TV(\mathbf{u})
\end{equation}
Thus, the combinatorial version of problem \eqref{eq:tvlalpha} on an arbitrary lattice reads
\begin{equation}\label{eq:tvlalphadiscrete}\tag{TV$\ell^\alpha$}
		\min_{\mathbf{u}\in\mathcal{L}^N} F_\alpha(\mathbf{u}).
\end{equation}
This is the central problem of this paper. As already pointed out in the preceding section, functional $F_\alpha$ has the desirable property of being submodular, that is its pairwise terms $g(u_p,u_q) = \omega_{pq}\left|u_p - u_q\right|$ satisfy the inequality
\begin{equation}\label{eq:submod}
		g(x\lor y) + g(x\land y) \leq g(x) + g(y), \quad x,y \in \mathcal{L}^2,
\end{equation}
where $\lor$, $\land$ denote the element-wise maximum and minimum, respectively \cite{Mur03}. That \eqref{eq:submod} indeed holds can be shown by using the discrete coarea formula \cite{ChaDar09}.

Submodularity implies that multi-label graph cut algorithms can be applied directly to solve problem \eqref{eq:tvlalphadiscrete} either approximately \cite{BoyVekZab01} or exactly \cite{Ish03,Dar09}. However, a discrete TV has (by definition) the even stronger property \eqref{eq:discretecoarea}, which can be used to make solving \eqref{eq:tvlalphadiscrete} significantly more efficient by decoupling it into binary energies. We refer to \cite{DarSig06a,ChaDar09,GolYin09} on how to do so.
\subsection{Relations between the Continuous and Discrete Problems} \label{sec:relations}
The continuous TV denoising problem \eqref{eq:tvlalpha} is an infinite-dimensional one. In order to be numerically solvable, it must be brought to finite dimensions. This is often done by approximating $BV(\Omega)$ with a sequence of finite-dimensional subspaces. The general theoretical framework for the discretization of variational regularization in a Banach space setting can be found in \cite{PoeResSch10}. In the following paragraphs, we briefly examine problem \eqref{eq:tvlalphadiscrete} within this context.

A simple subspace of $BV(\Omega)$ is $V(\Omega_B^d)$: the class of functions that are piecewise constant on the Voronoi cells of lattice $\Omega_B^d$. Assume for simplicity $u\in V(\Omega_B^d)$ to be quantized, i.e.
\begin{equation}
	u = \sum_p{u_p\chi^{V_p}},\quad u_p\in\mathcal{L},
\end{equation}
where $\chi^{V_p}$ is the characteristic function of the Voronoi cell of lattice point $p$. Using the coarea formula \eqref{eq:coarea} while observing that level contours of a piecewise constant function can be computed by adding side lengths of Voronoi cells, we obtain \cite{CasKunPol99}:
\begin{equation} \label{eq:piecewiseconstanttv}
	\TV(u) = \sum_{l=0}^{L-2} \sum_{p,q} \frac{1}{2}\abs{V_p \cap V_q} \abs{\chi^l_p - \chi^l_q} =  \frac{1}{2} \sum_{p,q} \abs{V_p \cap V_q} \abs{u_p - u_q}.
\end{equation}
The factor $\frac{1}{2}$ accounts for the fact that every pair of neighbouring cells appears twice in the sum.

Now consider the sampled version of $u\in V(\Omega_B^d)$, that is, a discrete image $\mathbf{u}$ defined on $\Omega_B^d$. Let either $B=I$ or $B=H$, so that the Voronoi cells are regular polygons, with side length $v$. Additionally, let the boundary of $\Omega$ coincide with boundaries of Voronoi cells\footnote{If this is not the case, the following observation holds up to a negligible constant accounting for boundary effects.} and denote by $\sim$ the nearest-neighbour relation, i.e.\ $\abs{V_p \cap V_q} = v$, iff $p \sim q$. Then, using \eqref{eq:piecewiseconstanttv}, the discrete functional \eqref{eq:discrete-functional} can be rewritten as
\begin{align} \label{eq:connection}
	F_\alpha(\mathbf{u})
	& = \lambda \sum_{p}\abs{V_p}\abs{u_p-f_p}^\alpha + \omega \sum_{p \sim q} \abs{u_p - u_q} \notag \\
	& = \lambda \norm{u - \tilde{f}}^\alpha_{L^\alpha} + \frac{\omega}{v} \TV(u),
\end{align}
where $\tilde{f}$ is the piecewise constant extension of $\mathbf{f} = (f_p)_p$, in other words, the piecewise constant approximation of data $f$. Note that for the nearest-neighbour relation, the Cauchy-Crofton formula prescribes only \emph{one} weight $\omega$. We conclude that, under the above presumptions, minimization of $F_\alpha$ can be interpreted as choosing $V(\Omega_B^d)$ as approximation space to the continuous problem, that is, solving $\eqref{eq:tvlalpha}$ restricted to $V(\Omega_B^d) \subset BV(\Omega)$ only with a different regularization parameter $\tilde{\lambda} = \frac{v}{\omega}\lambda$. Therefore, approximation properties for the space $V(\Omega_B^d)$ are directly passed down to the corresponding discrete problem.

In \cite{FitKee97,CasKunPol99} such approximation properties for piecewise constant functions on rectangular patches (including $V(\Omega_I^d)$) were derived. It was shown that an arbitrary $u\in BV(\Omega)$ can be approximated by those simple functions as $d\rightarrow 0$, only if the $\ell^1$ norm on $\mathbb{R}^2$ is used to define the TV (cf.\ \eqref{eq:tv}), making it anisotropic.\footnote{In \cite{BelLus03} these results were extended to general $\ell^k$ norms, $k \in [1,\infty]$, yet at the expense of having to resort either to irregular triangulations of $\Omega$ or to piecewise linear functions.} This can be regarded as a decisive factor behind widely observed metrication errors. Roughly speaking, the shapes of the constant patches dictate the type of (anisotropic) TV which is eventually approximated. This, in turn, causes artefacts in the denoised image, since edges that are aligned with grid directions will be ``cheaper" than others. Similar approximation properties of TV regularization on a hexagonal lattice have not been derived yet, but are under investigation in \cite{PoeResSch12_report}.

If larger neighbourhoods are considered, the connection to piecewise constant functions does no longer hold. In this case the Cauchy-Crofton weights can be expected to reduce the above mentioned discretization effects for any underlying lattice (cf.\ Sec.\ \ref{sec:cutmetric}).\footnote{Another possibility of reducing metrication errors consists in using higher order pairwise terms $\abs{u_p-u_q}^\beta$, $\beta>1$ \cite{SinGraSinVid11}. This, however, abandons the TV approach to image restoration and is, therefore, out of the scope of this paper.} However, they still display one important property even for nearest-neighbour grids, which turns out to be useful for experimental issues. An easy calculation shows, that the ratio $\frac{\omega}{v}$ from \eqref{eq:connection} is equal to $\frac{\pi}{4}$ for both lattices $\Lambda(dI)$ and $\Lambda(d'H)$ and any choices of $d,d'>0$. This implies that, if the discrete problem is solved on both grids with the same $\lambda$, then the regularization parameters of the corresponding continuous functionals will be equal as well ($\tilde{\lambda} = \lambda\frac{4}{\pi}$). In other words, the Cauchy-Crofton weights ensure that the restored images obtained from the two different spatial discretizations correspond to the same degree of regularization in the continuous model, and thus can be directly compared (cf.\ Sec.\ \ref{sec:experiments}).
\section{Experiments}\label{sec:experiments}
Every single one of the following five experiments complies with a fixed procedure, which is explained below.
\begin{enumerate}
	\item We generate two 8-bit ($L=2^8$) ground truth images of the same resolution, one standard and one hexagonally sampled. For a synthetic test image this can be achieved by simply sampling it to two lattices of the same density (cf.\ Sec.~\ref{sec:lattices}). A natural test image has to be resampled to a hexagonal grid, since we do not have access to a device producing hexagonally sampled data.
	\item Both ground truth images are contaminated with the same type of noise.
	\item For fixed $\alpha$ (depending on the noise) and a fixed range of reasonable values for $\lambda$ (that has been determined before) both resulting versions of problem \eqref{eq:tvlalphadiscrete} are solved exactly with the sequential algorithm from \cite{DarSig06a} in combination with the max-flow algorithm from \cite{BoyKol04}.
	\item Finally, steps 2 and 3 are repeated several times to make up for different realizations of noise.
\end{enumerate}

Special emphasis lies on comparing TV regularization with respect to grid topology $\mathcal{N}_H^6$ on the one side, to that with respect to $\mathcal{N}_I^4$ and $\mathcal{N}_I^8$ on the other. Apart from visually comparing our results, the denoising quality will be measured using two different performance figures: first, the $\ell^1$ distance between ground truth $\mathbf{u}^\star$ and restored image $\hat{\mathbf{u}}$ divided by the number of lattice points (to make up for the fact that the square and hexagonal images will not have exactly the same amount of pixels)
\begin{equation} \label{eq:pf1}
	\frac{1}{N}\norm{\mathbf{u}^\star - \hat{\mathbf{u}}}_1,
\end{equation}
which can be interpreted as the mean absolute error per sampling point; and second, the ratio of correctly restored pixels
\begin{equation} \label{eq:pf2}
	\frac{1}{N} \Abs{\{p\in\Omega_B^d: u_p^\star = \hat{u}_p\}}.
\end{equation}
Since the second figure is a less reliable measure of denoising quality, e.g.\ in cases where the restored image suffers from a loss of contrast, it will be only regarded as an additional source of information, whereas more emphasis will be placed upon \eqref{eq:pf1} and visual appearance. At the end of this section the experiments' outcome is briefly discussed.
\paragraph{Experiment No.\ 1}
The contrast enhanced Shepp-Logan phantom (Fig.~\ref{fig:shepp-logan}) was chosen as test image, and sampled to a square image of size $256\times 256$ and to a hexagonal one of roughly the same resolution ($275\times 238$). After adding $60\%$ salt \& pepper noise, the images were denoised with an $L^1$ data term ($\alpha = 1$) and using the three topologies $\mathcal{N}_H^6$, $\mathcal{N}_I^4$ and $\mathcal{N}_I^8$ and the corresponding weights from Sec.~\ref{sec:cutmetric}. The procedure was repeated 50 times for each of several different values of $\lambda$. This setup, a piecewise constant image degraded with impulse noise is perfectly suited for TV$L^1$ restoration and can be expected to yield good results.

\begin{figure}[!ht]
  \centering
  \includegraphics[height=0.31\textwidth]{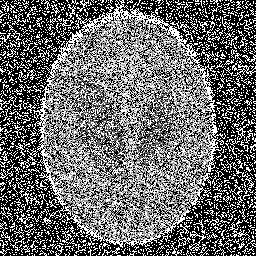}
  \hspace{1pt}
  \includegraphics[height=0.31\textwidth]{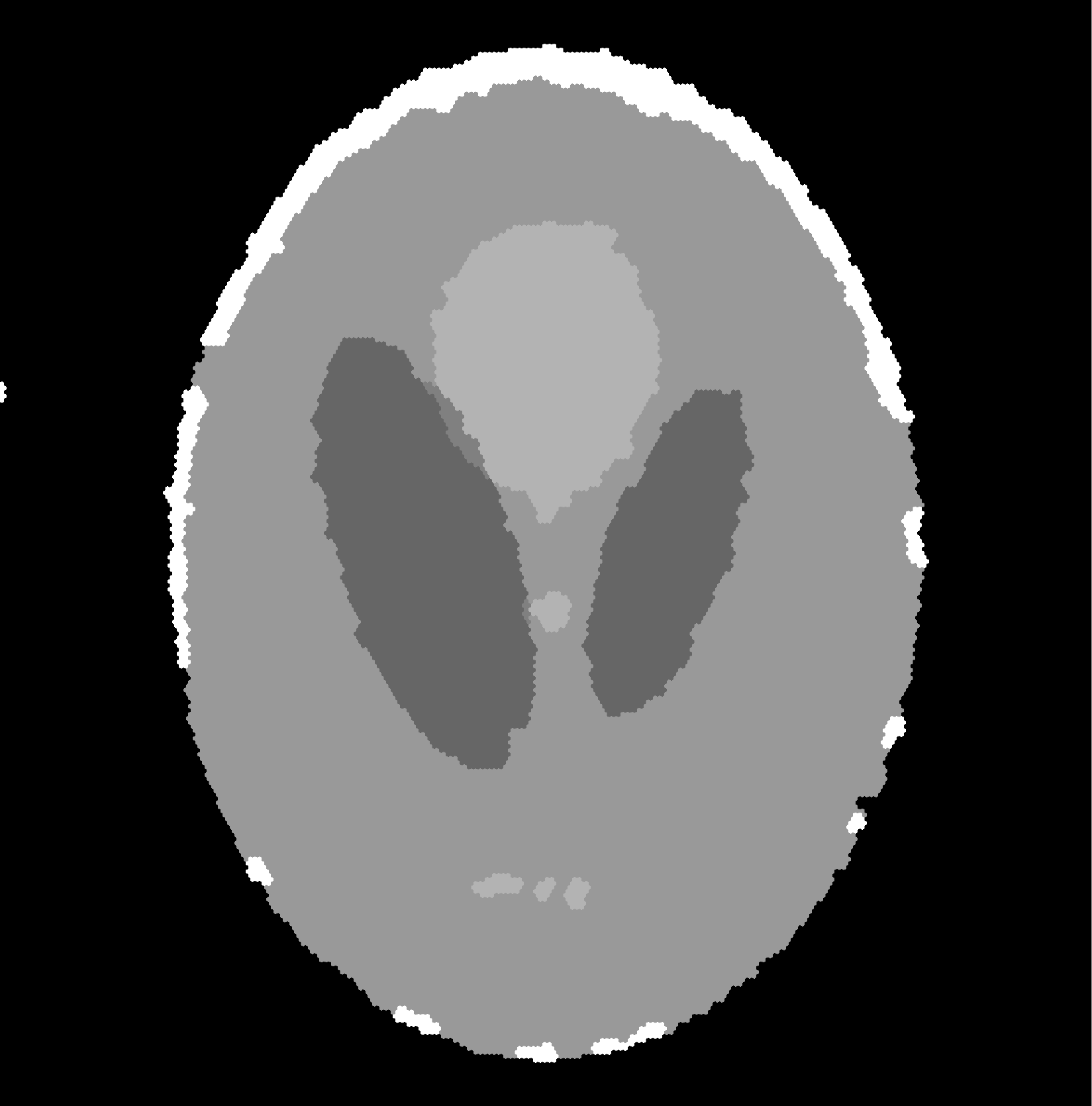}
  \hspace{1pt}
  \includegraphics[height=0.31\textwidth]{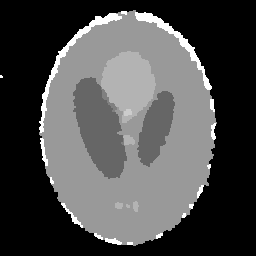}
  \caption{Results of Experiment No.~1. Noisy image opposed to denoised images for different grid topologies: $\mathcal{N}_H^6$ (centre), $\mathcal{N}_I^4$ (right). The values of $\lambda$ have been chosen to match with optimal performance in terms of $\ell^1$ distance to ground truth (cf.\ Fig.~\ref{fig:exp1}).}
  \label{fig:exp1img}
\end{figure}

\begin{figure}[!ht]
  \includegraphics[height=.49\textwidth]{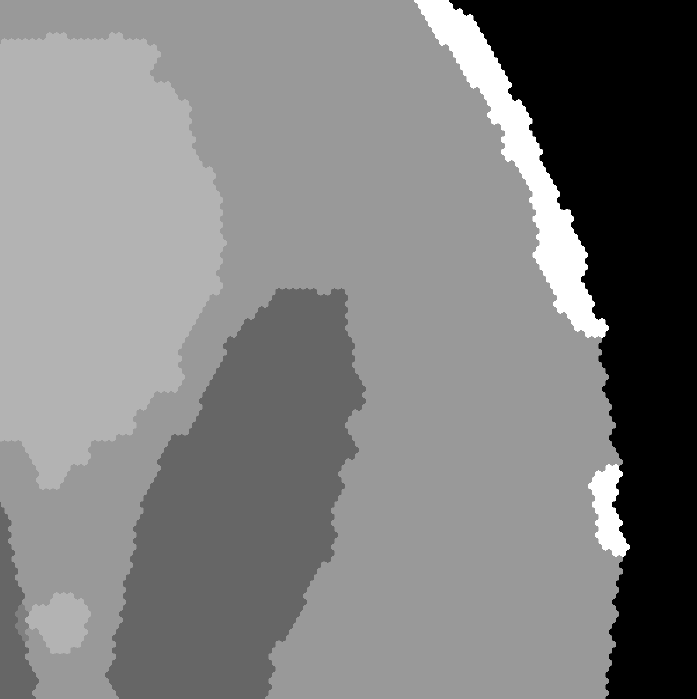}
  \hfill
  \includegraphics[height=.49\textwidth]{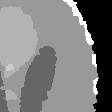}
  \caption{Results of Experiment No.~1. Details of restored images from Fig.~\ref{fig:exp1img}.}
  \label{fig:exp1detail}
\end{figure}

\begin{figure}[!ht]
  \includegraphics[height=0.21\textheight]{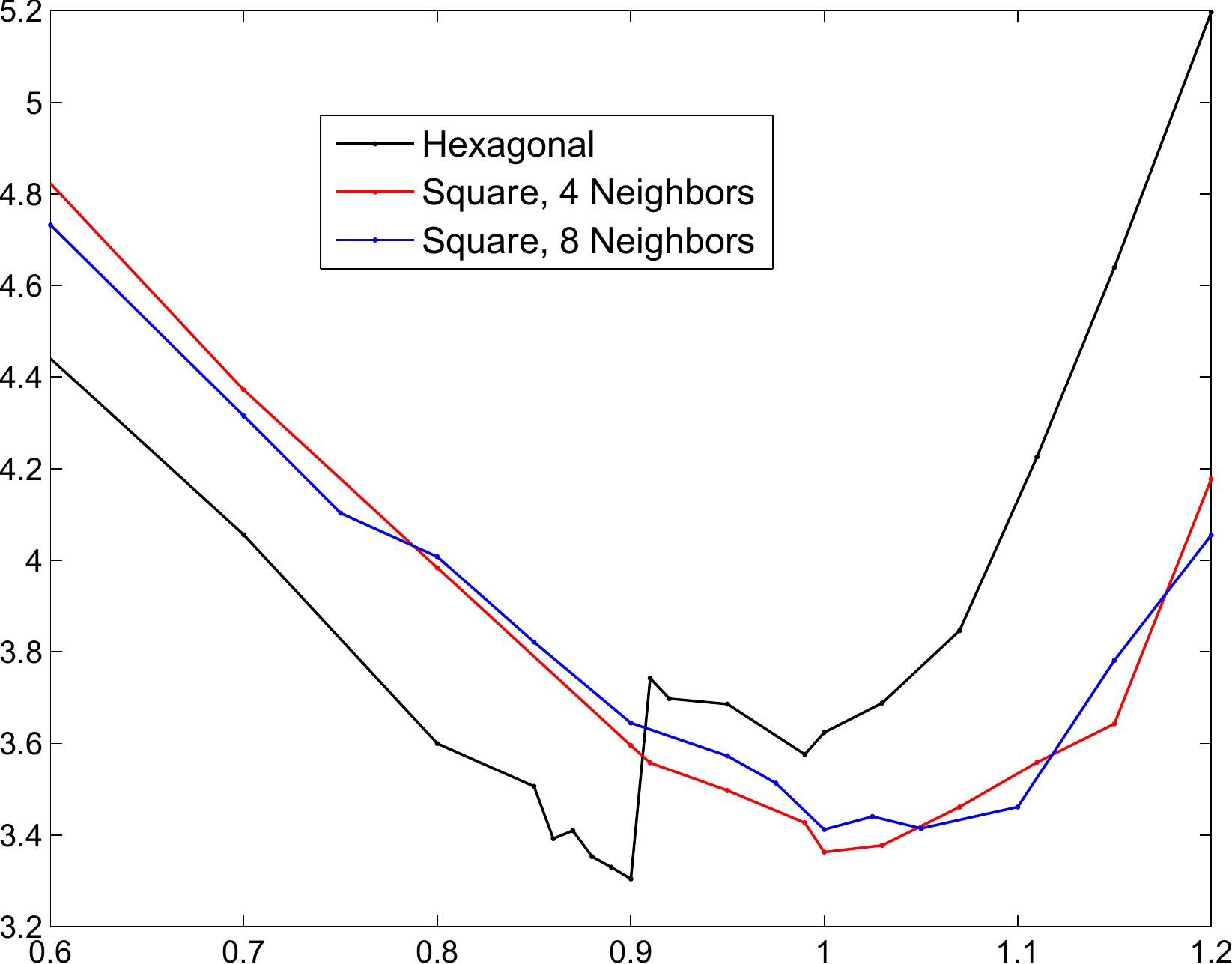}
  \hfill
  \includegraphics[height=0.21\textheight]{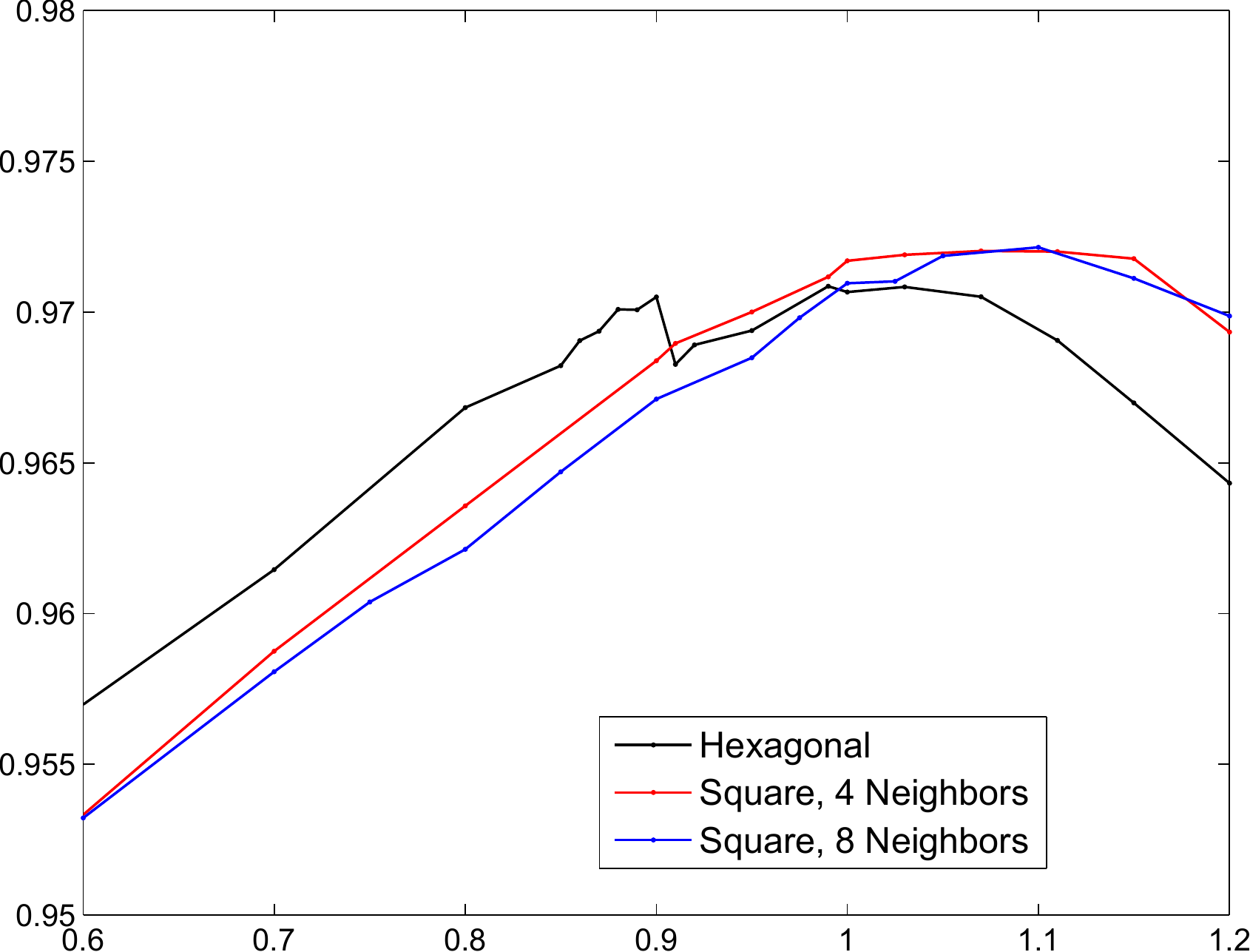}
  \caption{Results of Experiment No.~1. Averaged $\ell^1$ distances (left) and averaged ratios of correctly restored pixels (right) plotted against $\lambda$.}
  \label{fig:exp1}
\end{figure}

\begin{figure}[!ht]
  \includegraphics[width=0.49\textwidth]{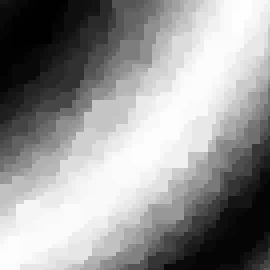}
  \hfill
  \includegraphics[width=0.49\textwidth]{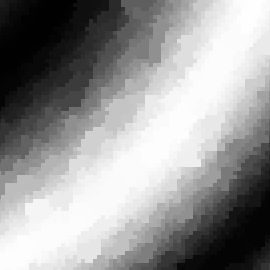}
  \\ \vskip1pt
  \includegraphics[width=0.49\textwidth]{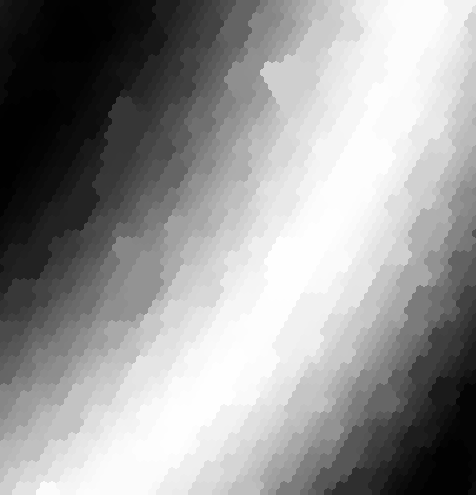}
  \hfill
  \includegraphics[width=0.49\textwidth]{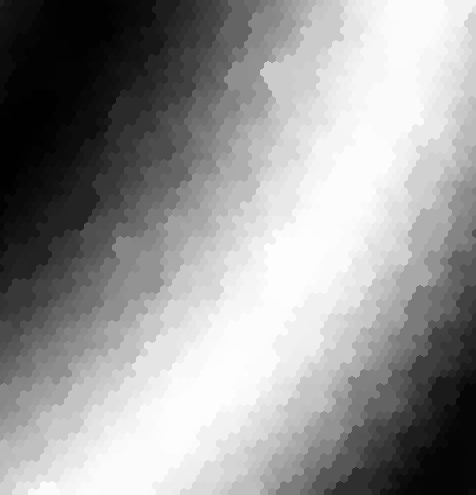}
  \caption{Results of Experiment No.~2. Magnified portions of denoised cosine test image (cf.\ Fig.~\ref{fig:cos}) for different topologies: $\mathcal{N}_I^{4}$ (top left), $\mathcal{N}_I^{16}$ (top right), $\mathcal{N}_H^6$ (bottom left), $\mathcal{N}_H^{12}$ (bottom right). Values of $\lambda$ have again been chosen to roughly match with optimal performance according to Fig.~\ref{fig:exp2}.}
  \label{fig:exp2img}
\end{figure}

\begin{figure}[!ht]
  \includegraphics[height=0.21\textheight]{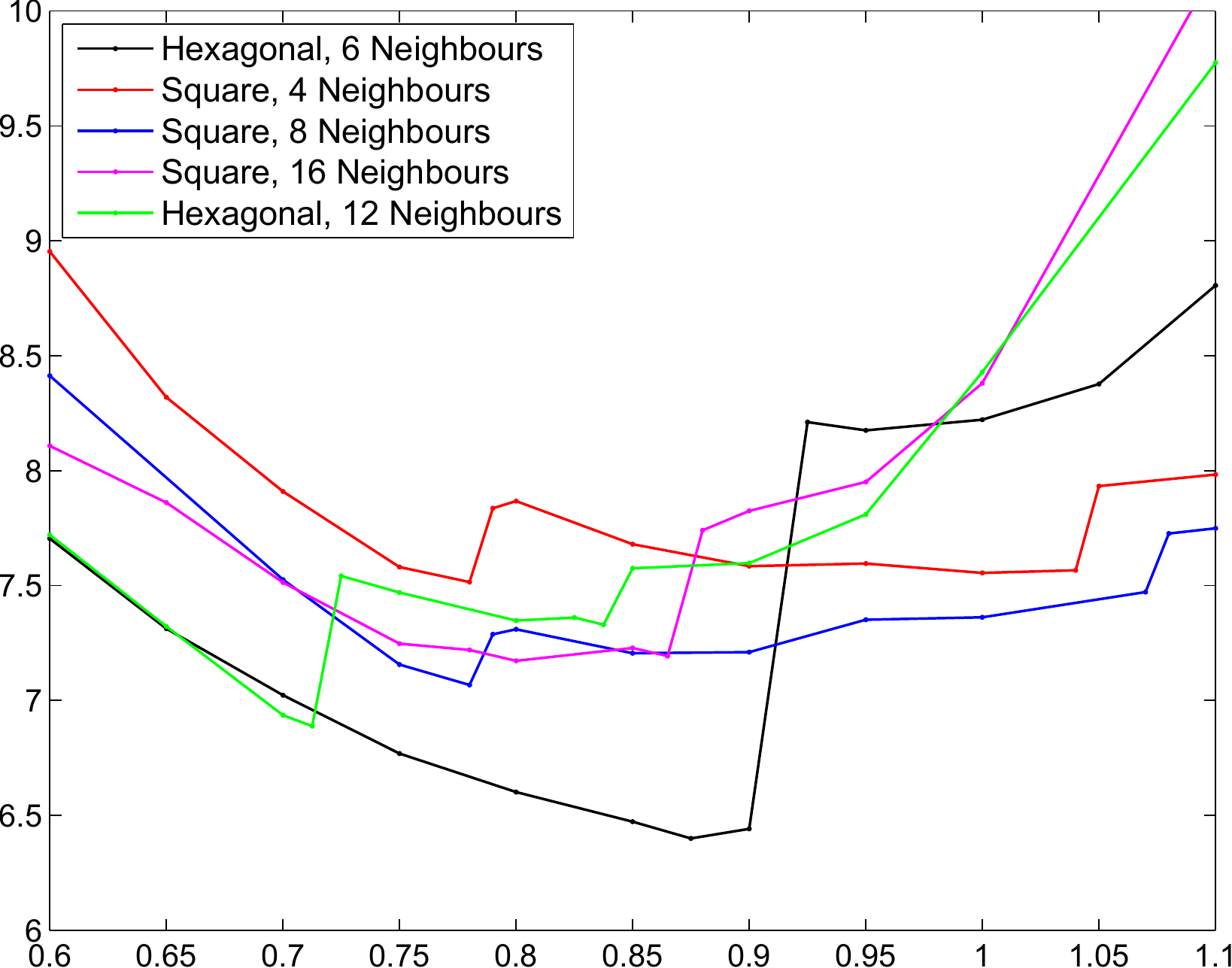}
  \hfill
  \includegraphics[height=0.21\textheight]{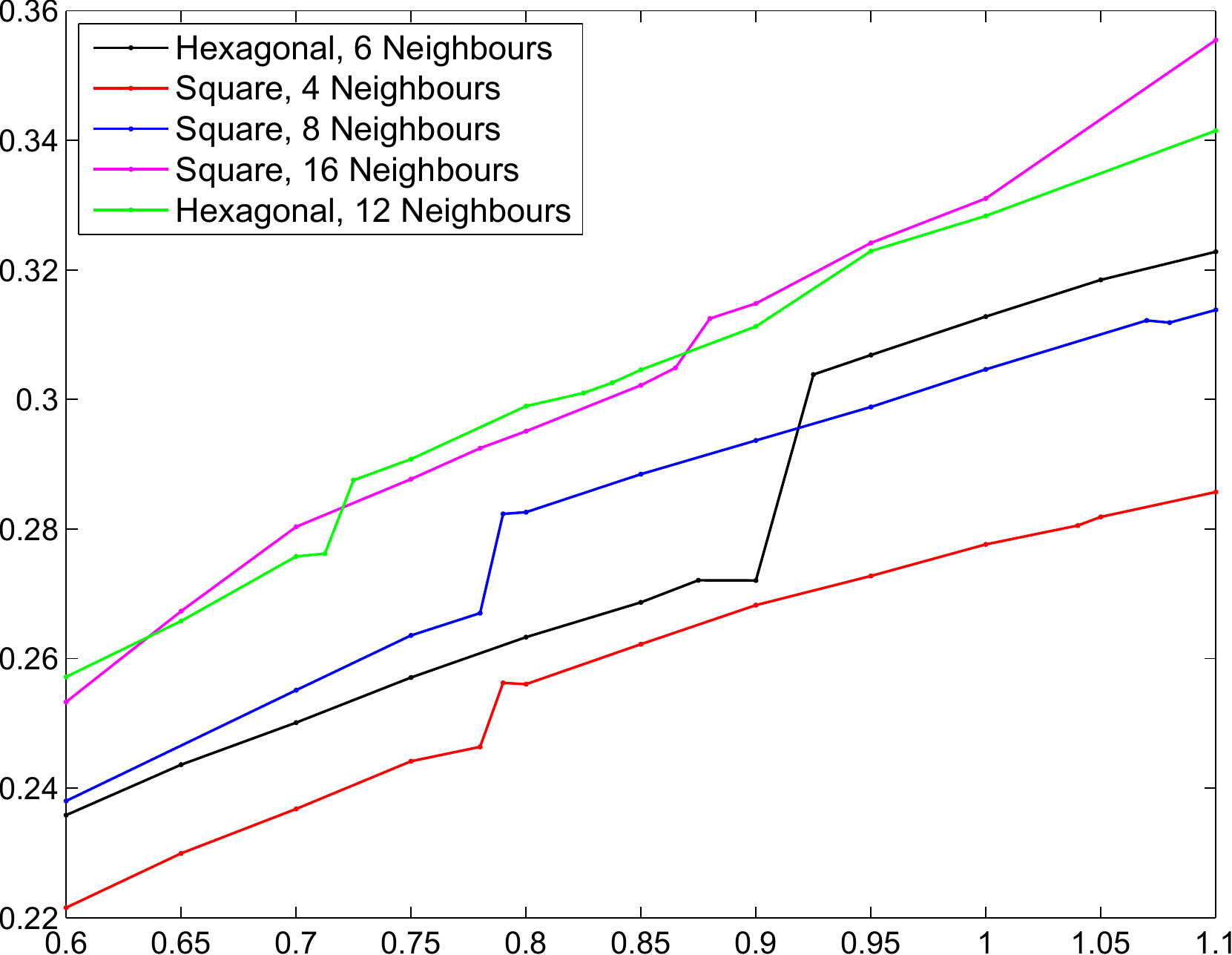}
  \caption{Results of Experiment No.~2. Averaged $\ell^1$ distances (left) and averaged ratios of correctly restored pixels (right) plotted against $\lambda$.}
  \label{fig:exp2}
\end{figure}

While the square grid topologies achieved optimal results in terms of \eqref{eq:pf1} and \eqref{eq:pf2} for $\lambda$ close to 1, for the hexagonal grid the best choice turned out to be $\lambda=0.9$ (cf.\ Fig.~\ref{fig:exp1}). This means that ``optimally" restored square lattice images are less filtered than optimal hexagonal images, since a higher regularization parameter means higher data fidelity and less smoothness (cf.\ the more ragged appearance of the square pixel image in Fig.~\ref{fig:exp1img} compared to the hexagonal one). If, however, a slightly smaller $\lambda$ is chosen for the square image, then certain image details (e.g.\ the tiny ellipses at the bottom) will disappear. Of course, this phenomenon depends on the specific test image and, to some extent, the raggedness is also caused by the more edgy, i.e.\ less circular, shape of square pixels. More importantly, Fig.~\ref{fig:exp1detail} indicates that TV regularization on a hexagonal lattice is less prone to metrication artefacts.
\paragraph{Experiment No.\ 2}
Experiment No.\ 1 was repeated with a different image: the cosine from Fig.~\ref{fig:cos} with a resolution of $270\times270$ pixels. The noise and $L^1$ data term were left unchanged. Since, in this experiment, the differences between $\mathcal{N}_I^{4}$, $\mathcal{N}_I^{8}$ and $\mathcal{N}_H^6$ were most pronounced (cf.\ Fig.~\ref{fig:exp2}), the larger neighbourhoods $\mathcal{N}_I^{16}$ and $\mathcal{N}_H^{12}$ were tested as well. Results are presented in Figs.~\ref{fig:exp2img} and \ref{fig:exp2}.
In contrast to the Shepp-Logan phantom the cosine test image consists entirely of smoothly varying intensity changes. TV restoration, however, is known to produce cartoon-like images, i.e.\ minimizers tend to be composed of subregions of more or less constant intensity separated by clear edges. This behaviour causes the so-called staircasing effect, which is clearly visible in Fig.~\ref{fig:exp2img} for both $\mathcal{N}_I^{4}$, where the restoration quality is additionally deteriorated by grave metrication errors, but also to some extent for $\mathcal{N}_H^6$. A comparison of $\mathcal{N}_I^{16}$ to $\mathcal{N}_H^{12}$ indicates that, by increasing neighbourhood sizes, those effects can be more effectively mitigated on the hexagonal lattice.

\begin{figure}[!ht]
  \centering
  \includegraphics[width=0.32\textwidth]{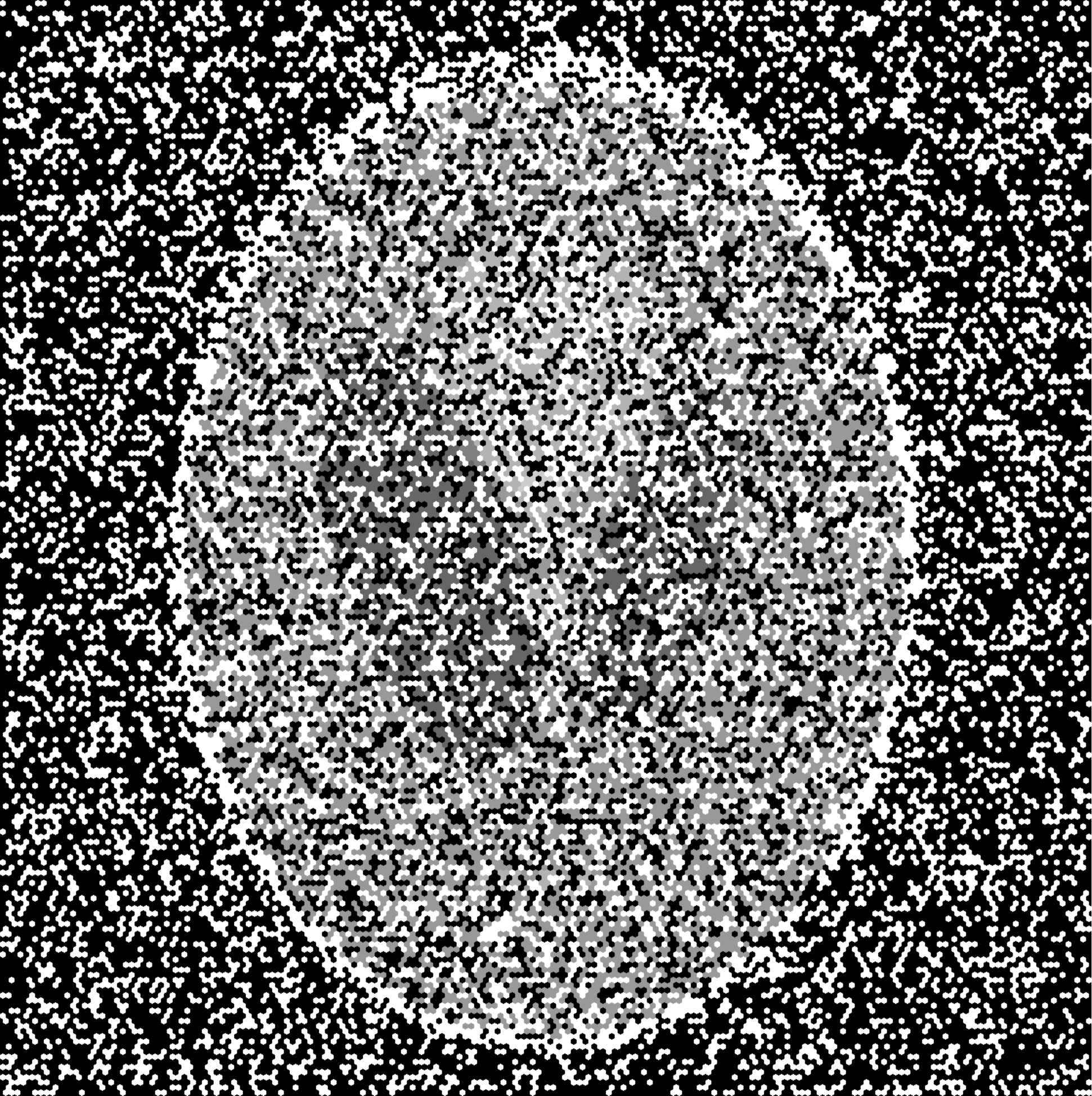}
  \hspace{1pt}
  \includegraphics[width=0.32\textwidth]{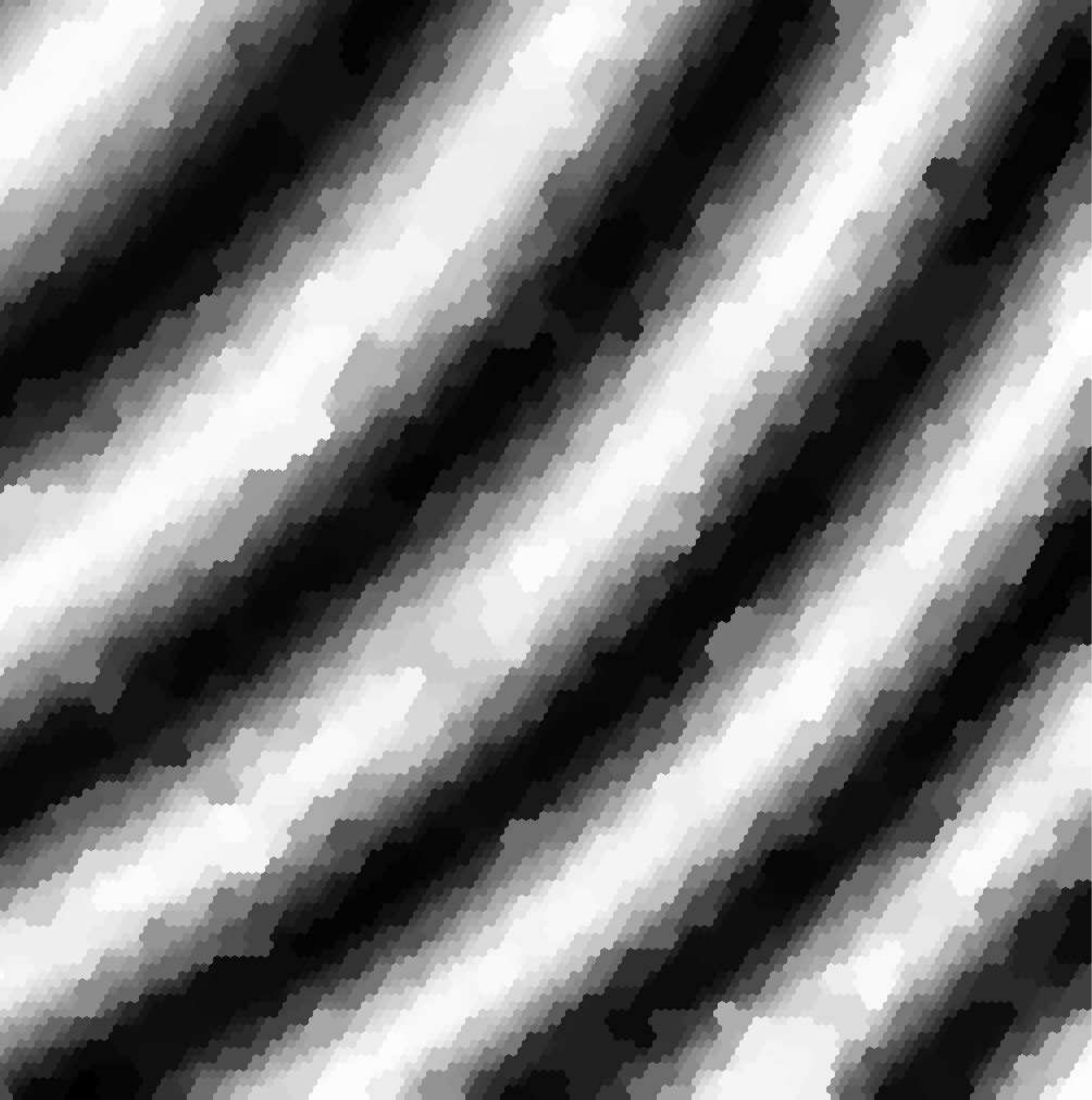}
  \hspace{1pt}
  \includegraphics[width=0.32\textwidth]{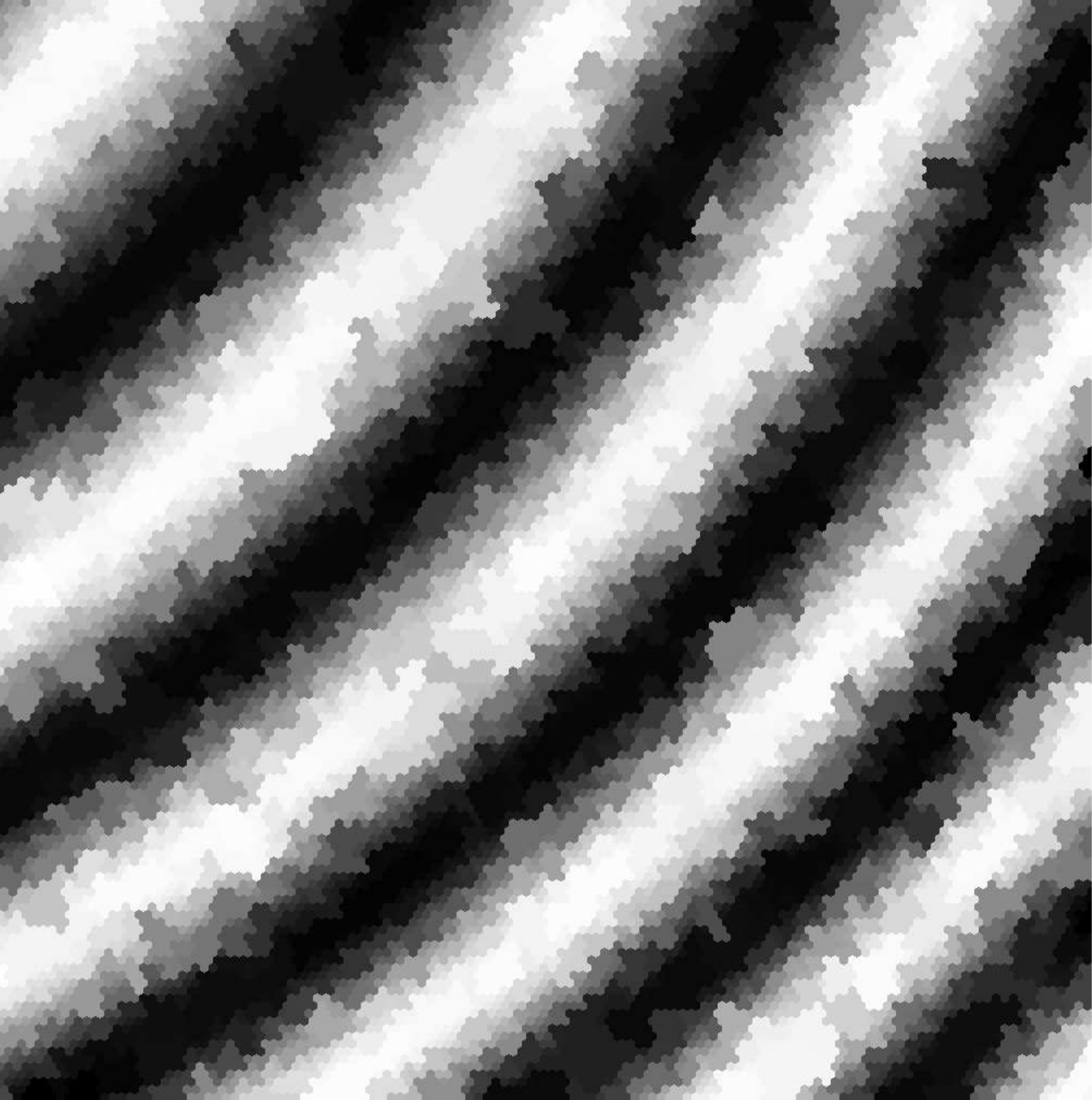}
  \caption{Results of Experiment No.~2. Hexagonally sampled cosine image corrupted with $60\%$ salt \& pepper noise (left), magnified portion of the same image denoised with regularization parameter $\lambda=0.906899$ (centre), and with $\lambda=0.9069$ (right). The image on the right, although of the same energy, is significantly less smooth; its TV is increased by about $24\%$, while its data fidelity is accordingly lower. This behaviour corresponds to the jump of the black line in Fig.~\ref{fig:exp2}.}
  \label{fig:jump}
\end{figure}

\begin{figure}[!ht]
  \includegraphics[width=0.49\textwidth]{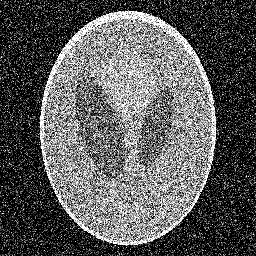}
  \hfill
  \includegraphics[width=0.49\textwidth]{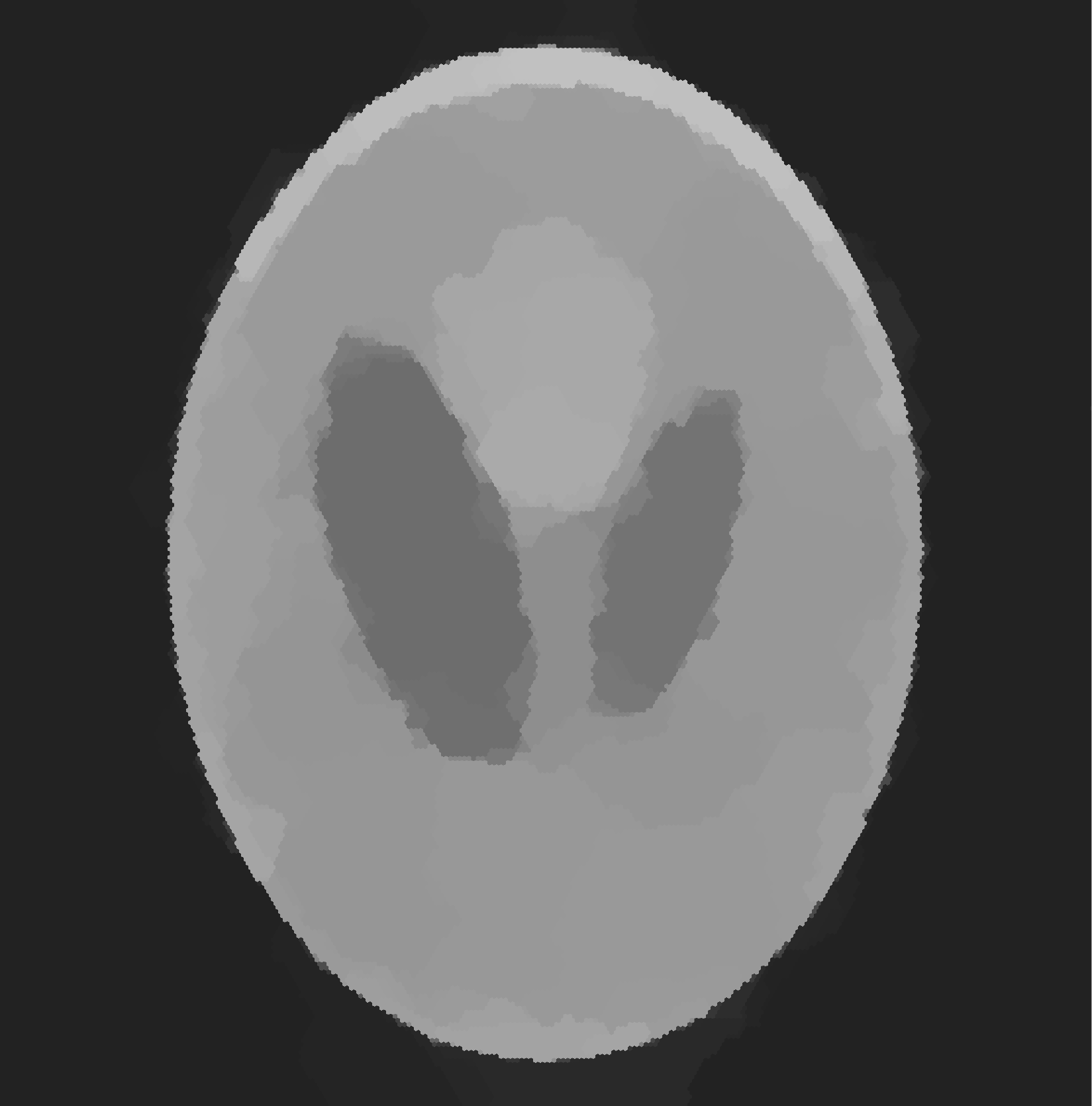}
  \\ \vskip1pt
  \includegraphics[width=0.49\textwidth]{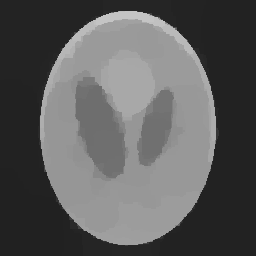}
  \hfill
  \includegraphics[width=0.49\textwidth]{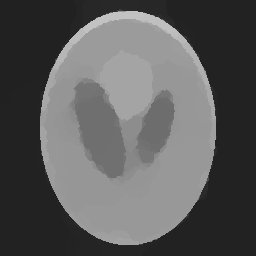}
  \caption{Results of Experiment No.~3. Noisy image (top left) opposed to denoised images with different grid topologies: $\mathcal{N}_H^{6}$ (top right), $\mathcal{N}_I^{4}$ (bottom left), $\mathcal{N}_I^8$ (bottom right). The regularization parameter $\lambda$ was set to $0.004$.}
  \label{fig:exp3img}
\end{figure}

\begin{figure}[!ht]
  \centering
  \includegraphics[height=0.21\textheight]{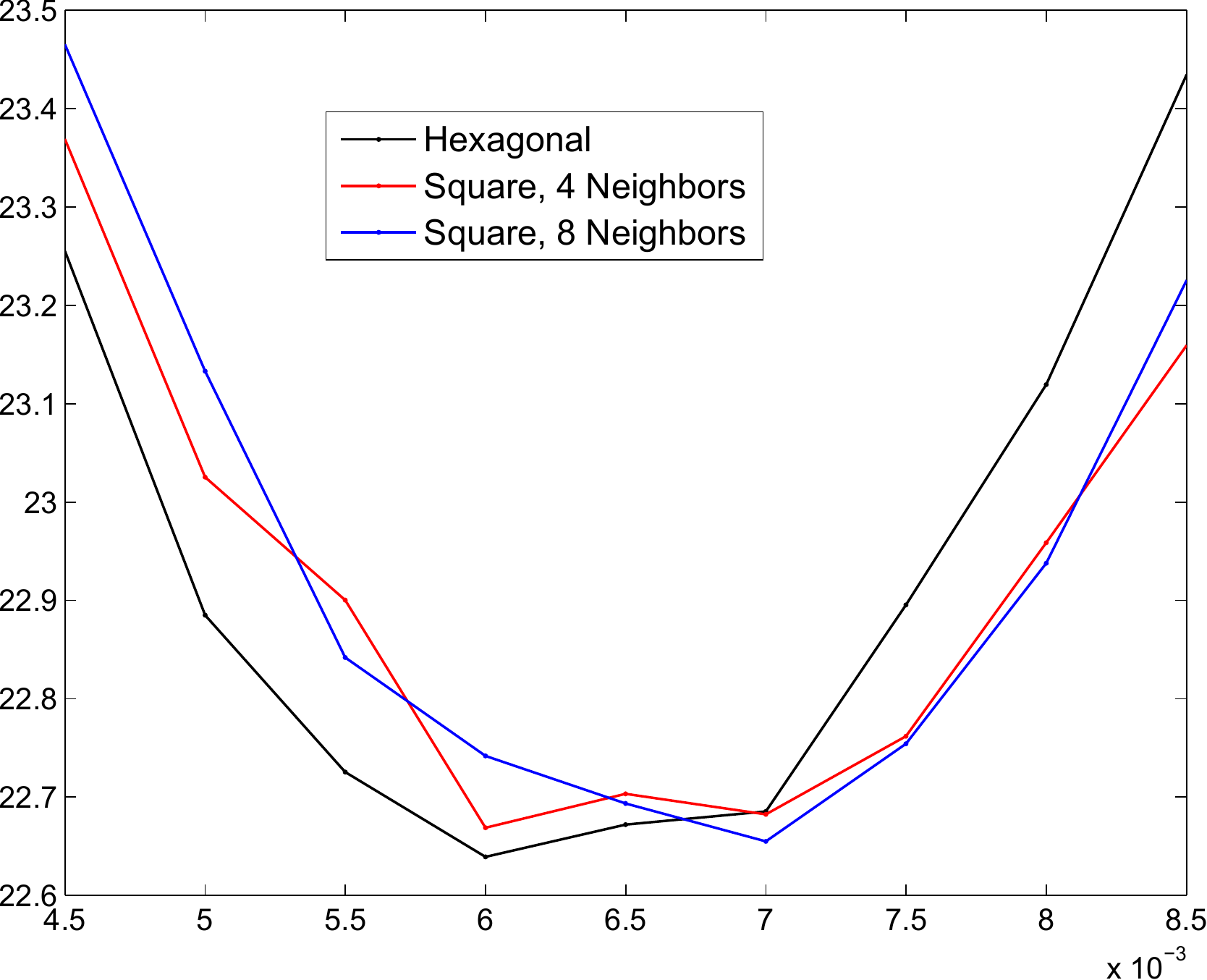}
  \caption{Experiment No.~3. Averaged $\ell^1$ distances plotted against $\lambda$.}
    \label{fig:exp3}
\end{figure}

It seems appropriate to make some additional remarks on this experiment in order to clarify its outcome and especially Fig.~\ref{fig:exp2}. The error curves in Fig.~\ref{fig:exp1} and even more so in Fig.~\ref{fig:exp2} display a striking feature: at certain distinguished values of the regularization parameter $\lambda$ they exhibit significant discontinuities. This peculiarity of TV$L^1$ regularization, has already been described in \cite{ChaEse05,SchGraGroHalLen09}. The authors showed that, in general, the data fidelity of minimizers depends discontinuously on $\lambda$, with at most countably many jumps. This behaviour, which is believed to be determined by the scales of image objects that rapidly merge at certain critical points, also manifests itself in other quantities, such as our performance figures. Fig.~\ref{fig:jump} elucidates this phenomenon.
\paragraph{Experiment No.\ 3}
The first experiment was repeated with the phantom being contaminated by additive Gaussian noise of zero mean and $10\%$ variance. Accordingly, an $L^2$ data term was employed to remove it. Since the TV$L^2$ model does not preserve contrasts, the number of correctly restored pixel intensities is close to zero for reasonable values of $\lambda$. We therefore only plot the $\ell^1$ distance in Fig.~\ref{fig:exp3}. This time the differences between grid topologies are less pronounced than in the previous experiments. In Fig.~\ref{fig:exp3img} results are presented.
\paragraph{Experiment No.\ 4}
Experiment No.\ 1 was repeated once more, but this time with a natural image, a $256\times256$ version of the well-known camera man. As stated at the beginning of this section, we need to resample natural images to a hexagonal grid. A rather simple and straightforward method was adopted in order to do so.

\begin{figure}[!ht]
  \includegraphics[width=0.48\textwidth]{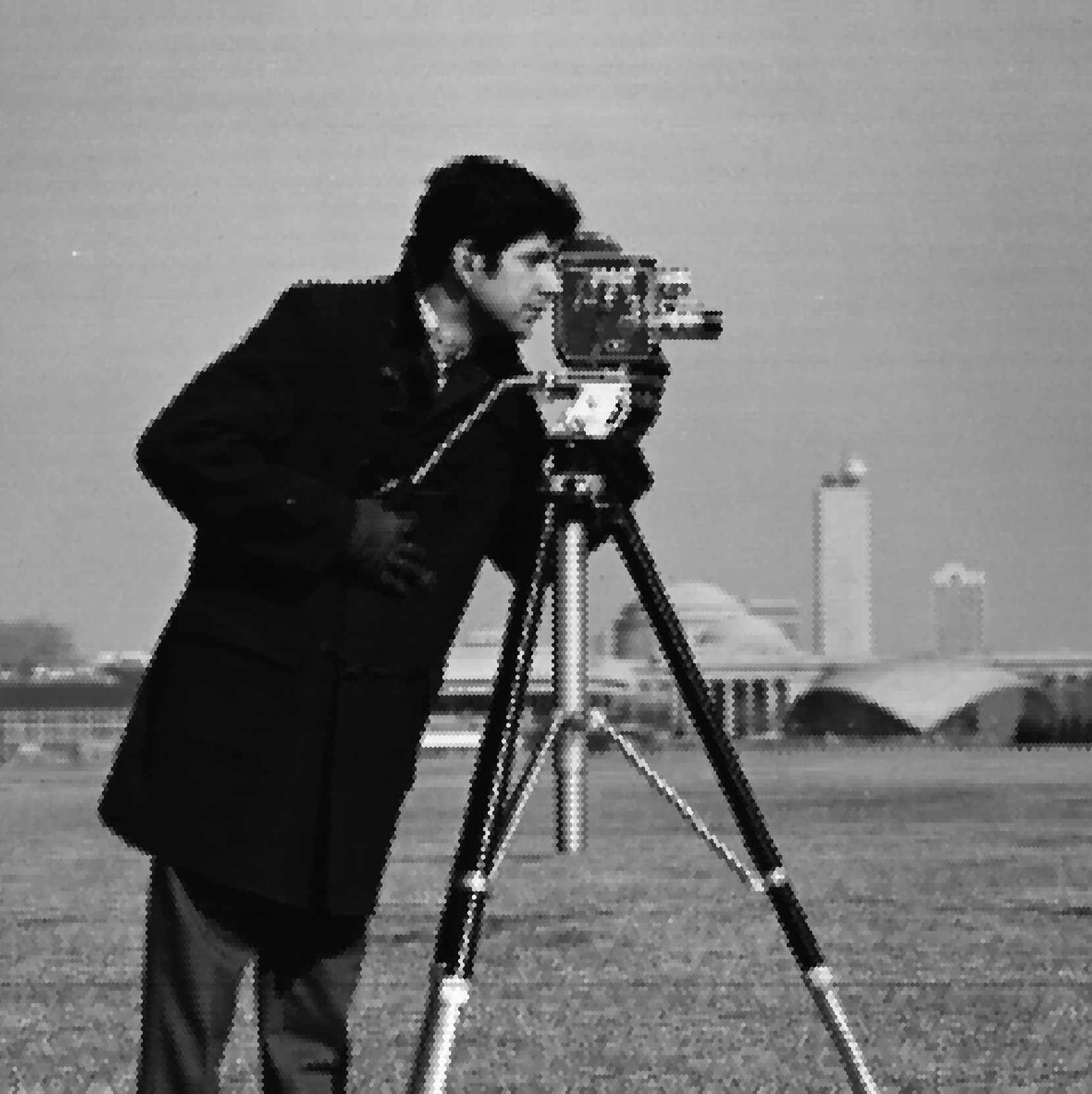}
  \hfill
  \includegraphics[width=0.48\textwidth]{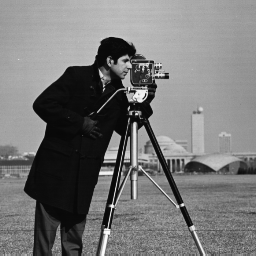}
  \caption{The two ground truth images used for Experiment No.\ 4. Original camera man image (right) and resampled to the hexagonal grid (left). It should be noted that the camera man features many edges that are perfectly aligned to Cartesian coordinates, e.g.\ contours of the camera, tripod and buildings in the background. Naturally those are better captured by a square pixel grid.}
  \label{fig:camera}
\end{figure}

First, the original image is enlarged by a certain scaling factor $c$, so that each pixel is replaced with a square of $c\times c$ pixels (if $c$ is integer). Now, the blown up image is covered with hexagonal hyperpixels, from which the pixel values of the hexagonal lattice are computed by taking some (weighted) average. If the hexagonal image thus generated is supposed to be of roughly the same resolution as the original one, the size of one hyperpixel should be close to $c^2$. For this experiment we chose $c=7.5$ and hyperpixels of size 56 \cite{MidSiv05}, which have also been used to create all the hexagonally sampled images in this note. In order to avoid smoothing effects as much as possible the median was employed instead of the arithmetic mean to compute the new pixel intensities. This results in a hexagonally sampled version of the camera man with a resolution of $274\times240$, having $0.34\%$ more sampling points than the original image (Fig.~\ref{fig:camera}).

\begin{figure}[!ht]
  \includegraphics[width=0.48\textwidth]{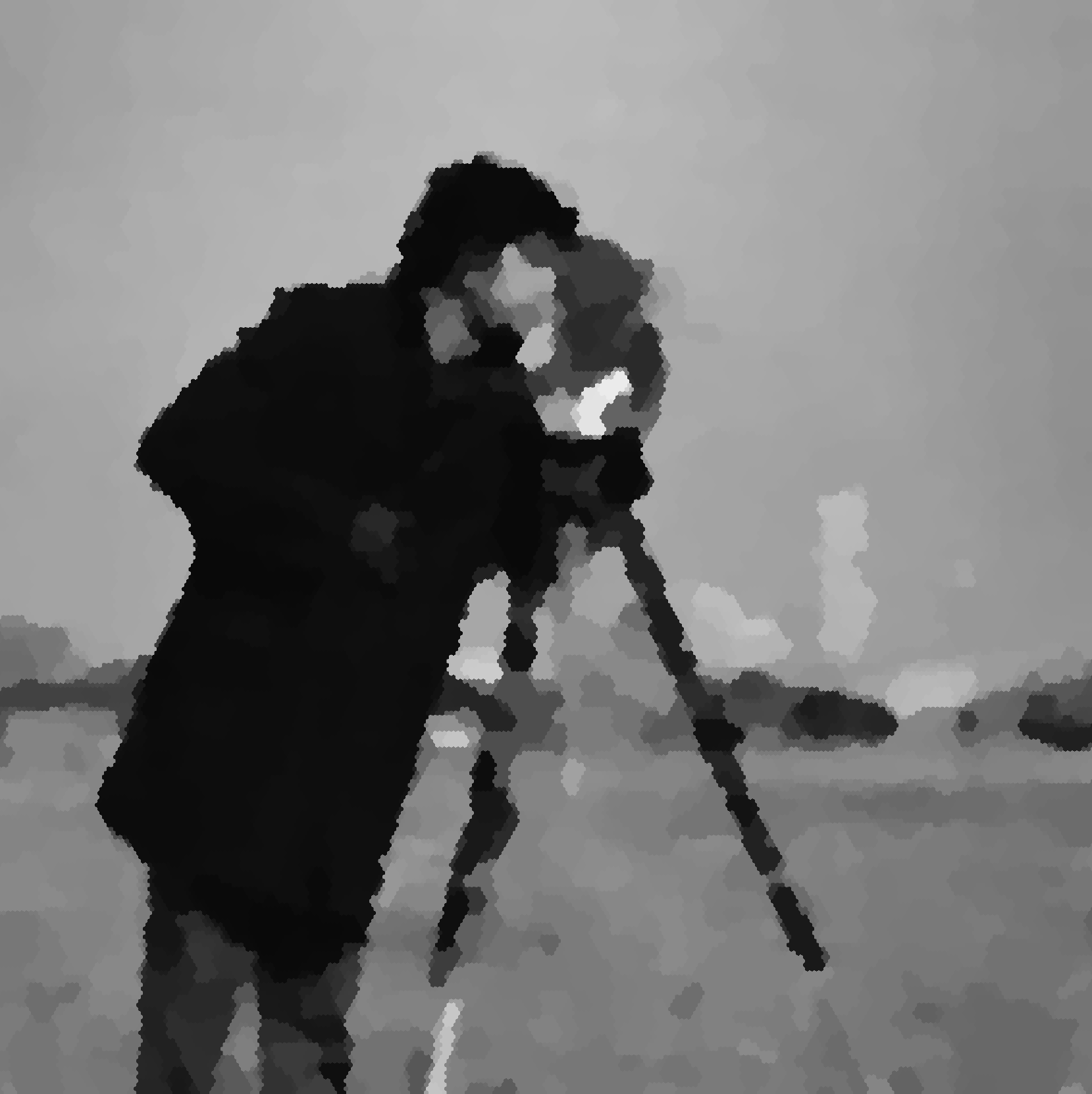}
  \hfill
  \includegraphics[width=0.48\textwidth]{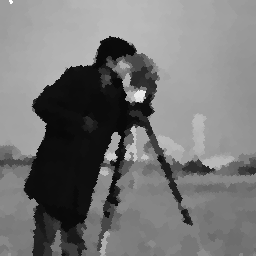}
  \caption{Results of Experiment No.~4. Denoised images with $\lambda=0.9$ and grid topologies $\mathcal{N}_H^6$ (left) and $\mathcal{N}_I^8$ (right).}
  \label{fig:exp4img}
\end{figure}

\begin{figure}[!ht]
  \includegraphics[height=0.21\textheight]{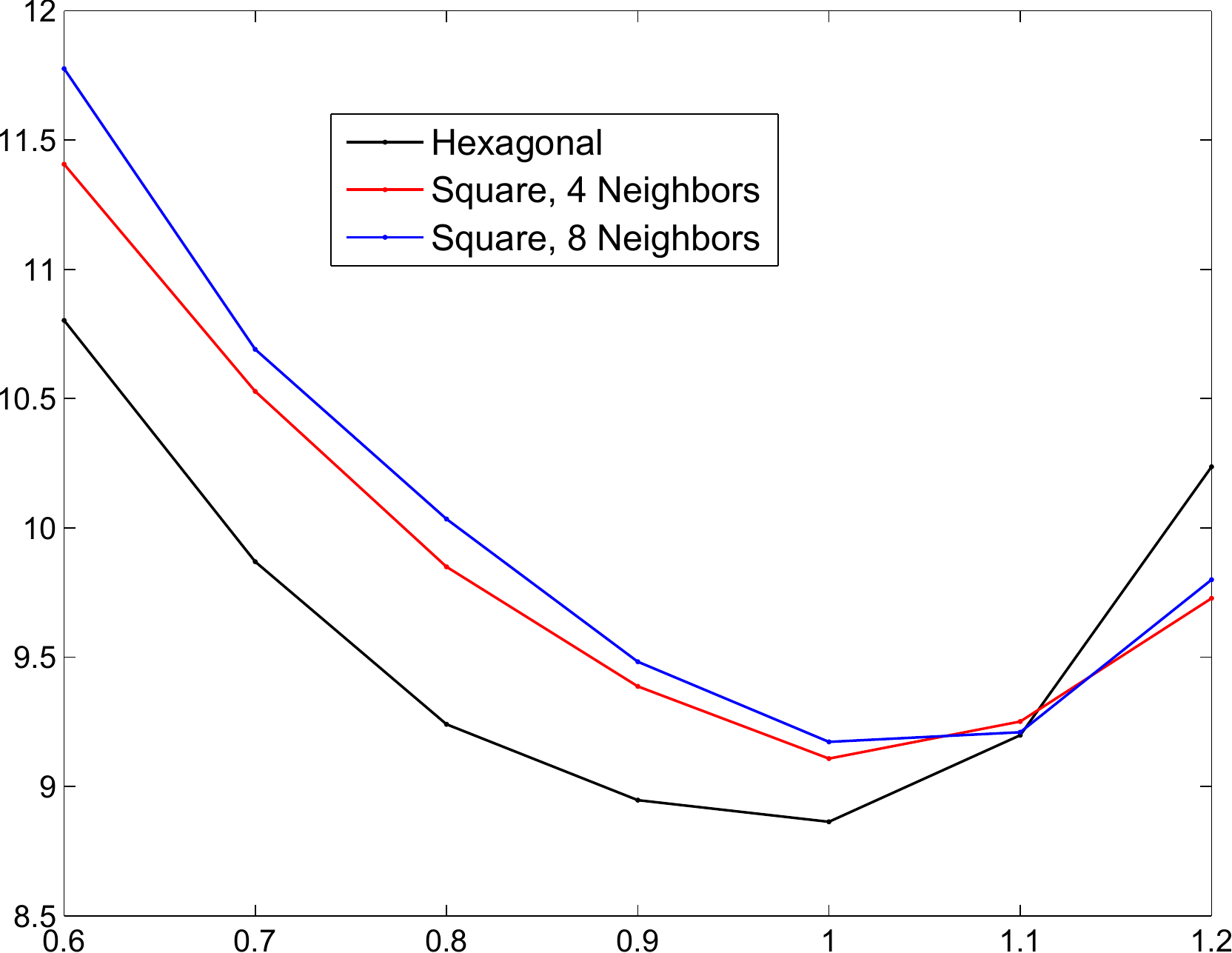}
  \hfill
  \includegraphics[height=0.21\textheight]{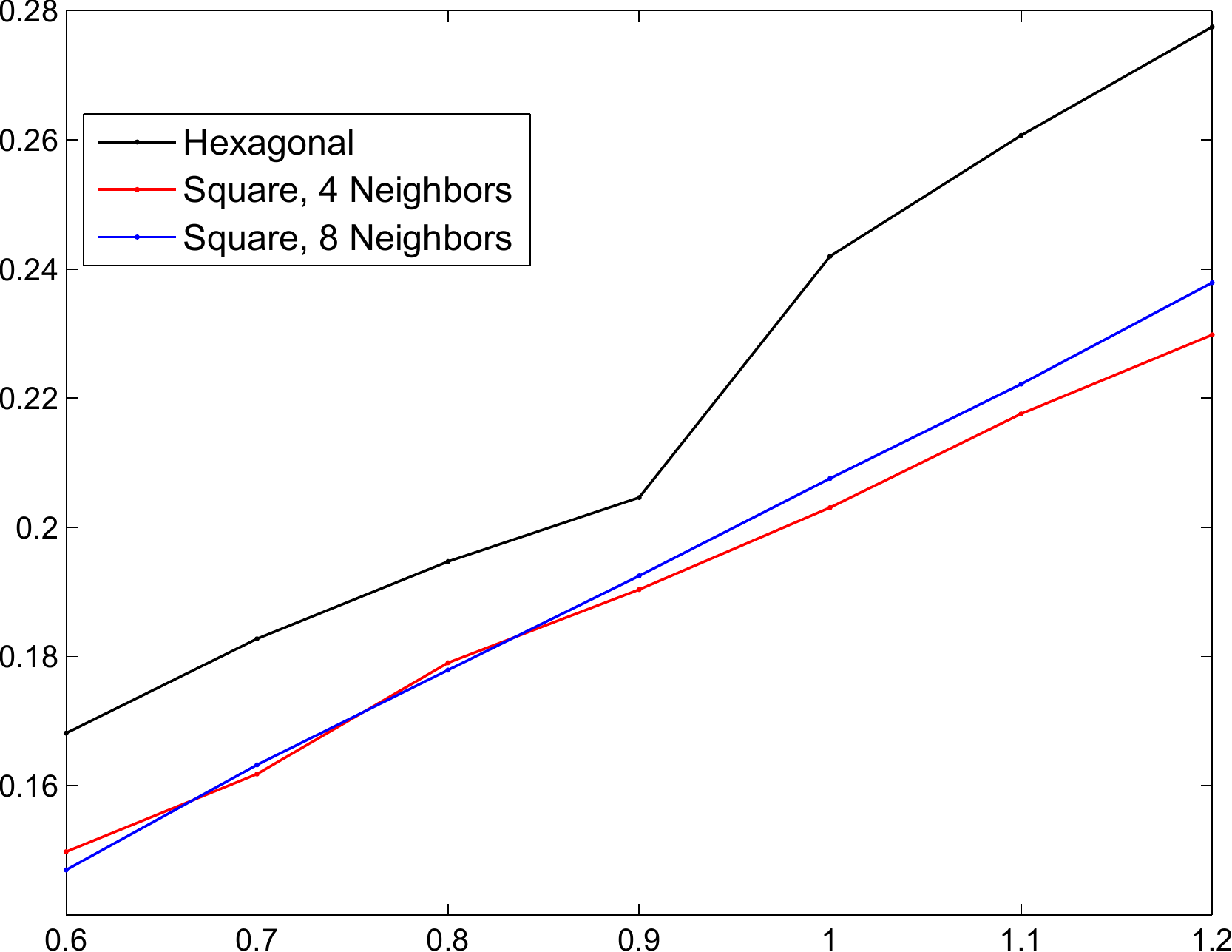}
  \caption{Results of Experiment No.~4. Averaged $\ell^1$ distances (left) and averaged ratios of correctly restored pixels (right) plotted against $\lambda$.}
  \label{fig:exp4}
\end{figure}

\begin{figure}[!ht]
  \includegraphics[width=0.49\textwidth]{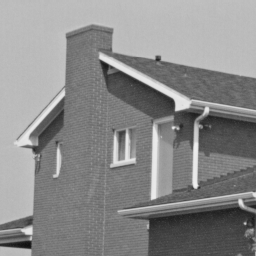}
  \hfill
  \includegraphics[width=0.49\textwidth]{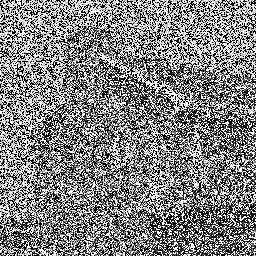}
  \\ \vskip1pt
  \includegraphics[width=0.49\textwidth]{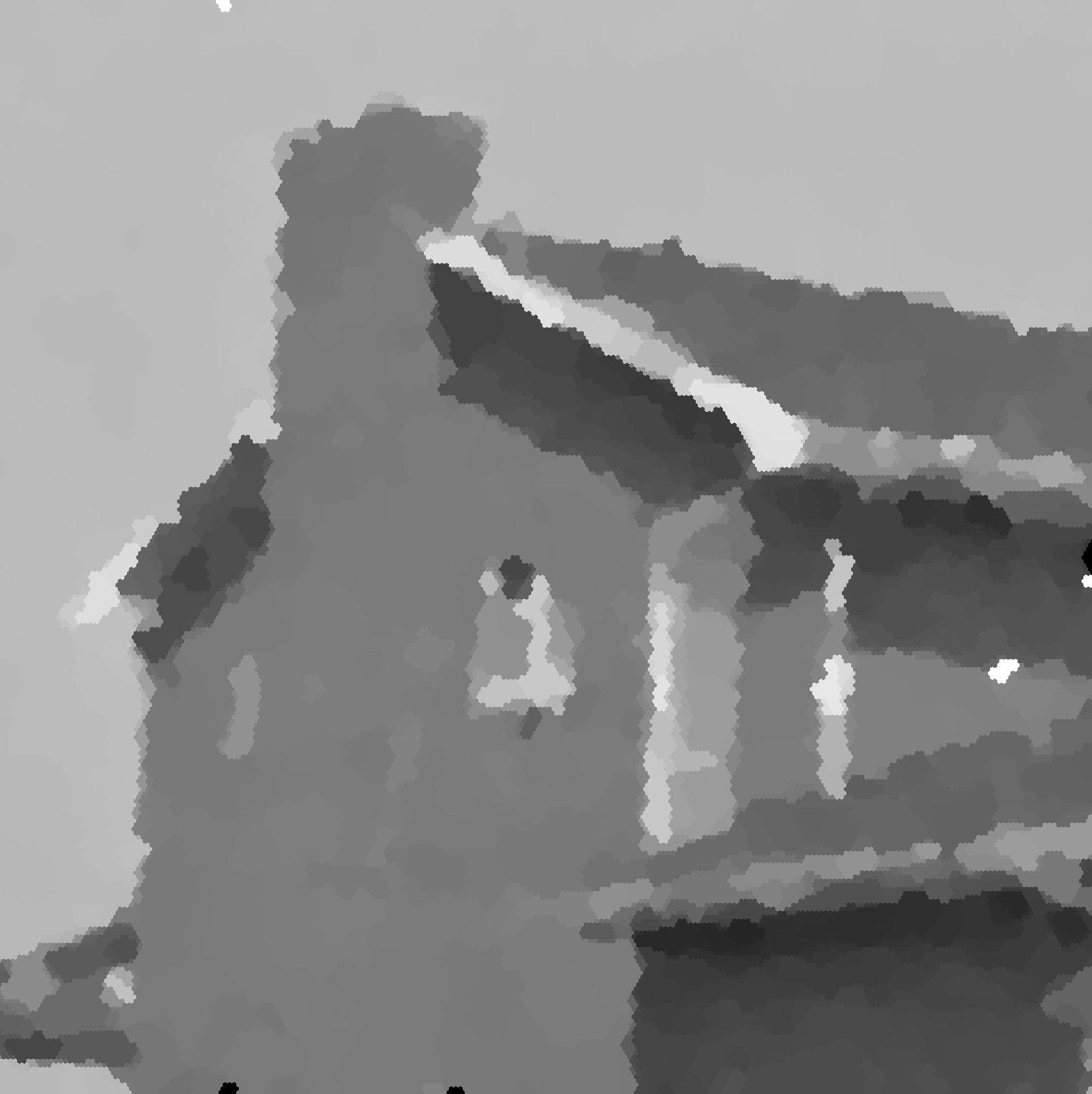}
  \hfill
  \includegraphics[width=0.49\textwidth]{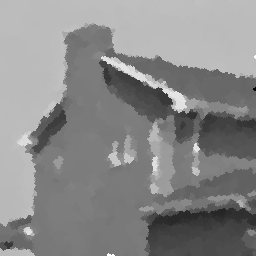}
  \caption{Results of Experiment No.~5. Square pixel ground truth image (top left) opposed to noisy image (top right) and restored images with respect to topologies $\mathcal{N}_H^6$ (bottom left) and $\mathcal{N}_I^8$ (bottom right). The regularization parameter $\lambda$ was set to 0.9.}
  \label{fig:exp5img}
\end{figure}

\begin{figure}[!ht]
  \includegraphics[height=0.21\textheight]{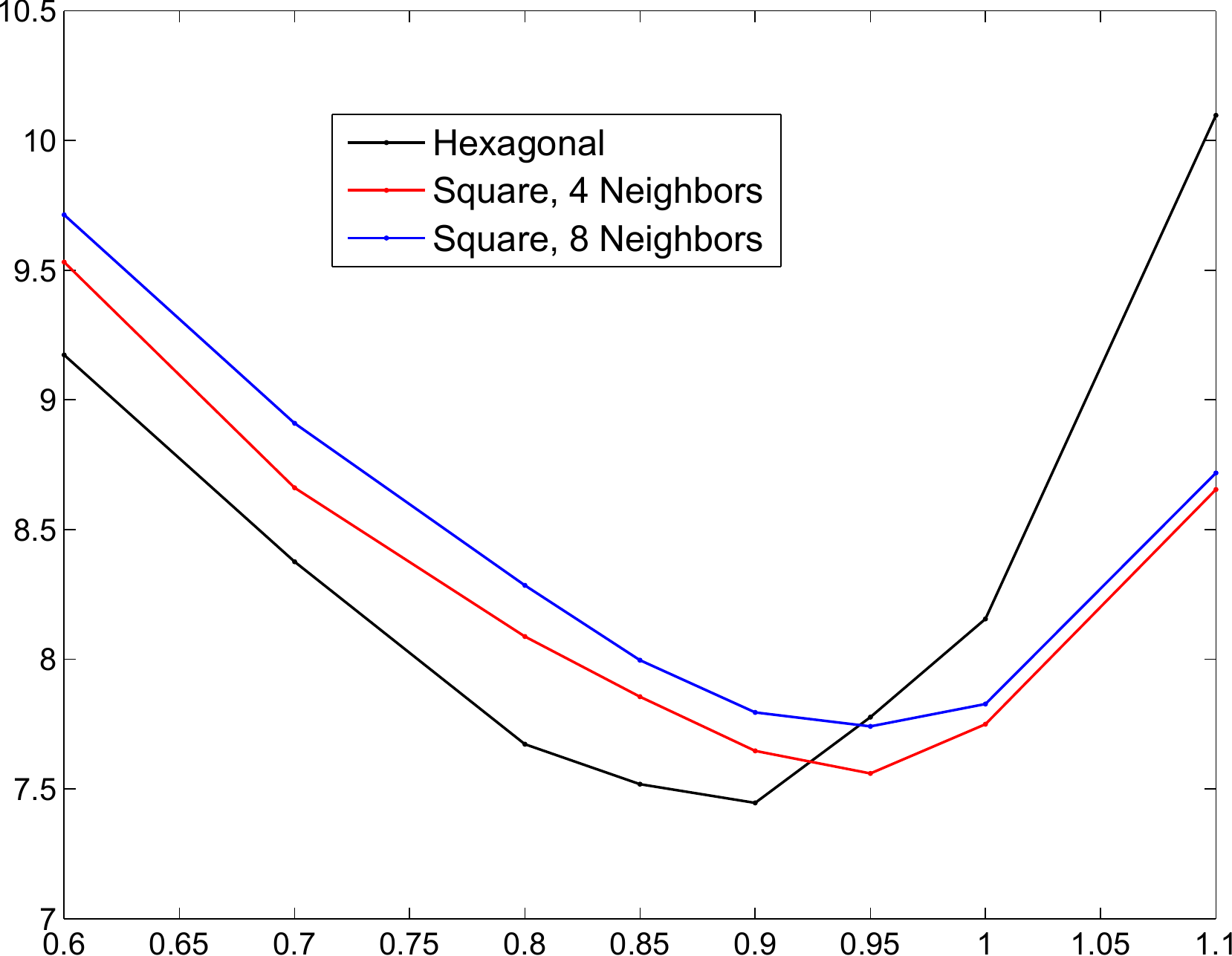}
  \hfill
  \includegraphics[height=0.21\textheight]{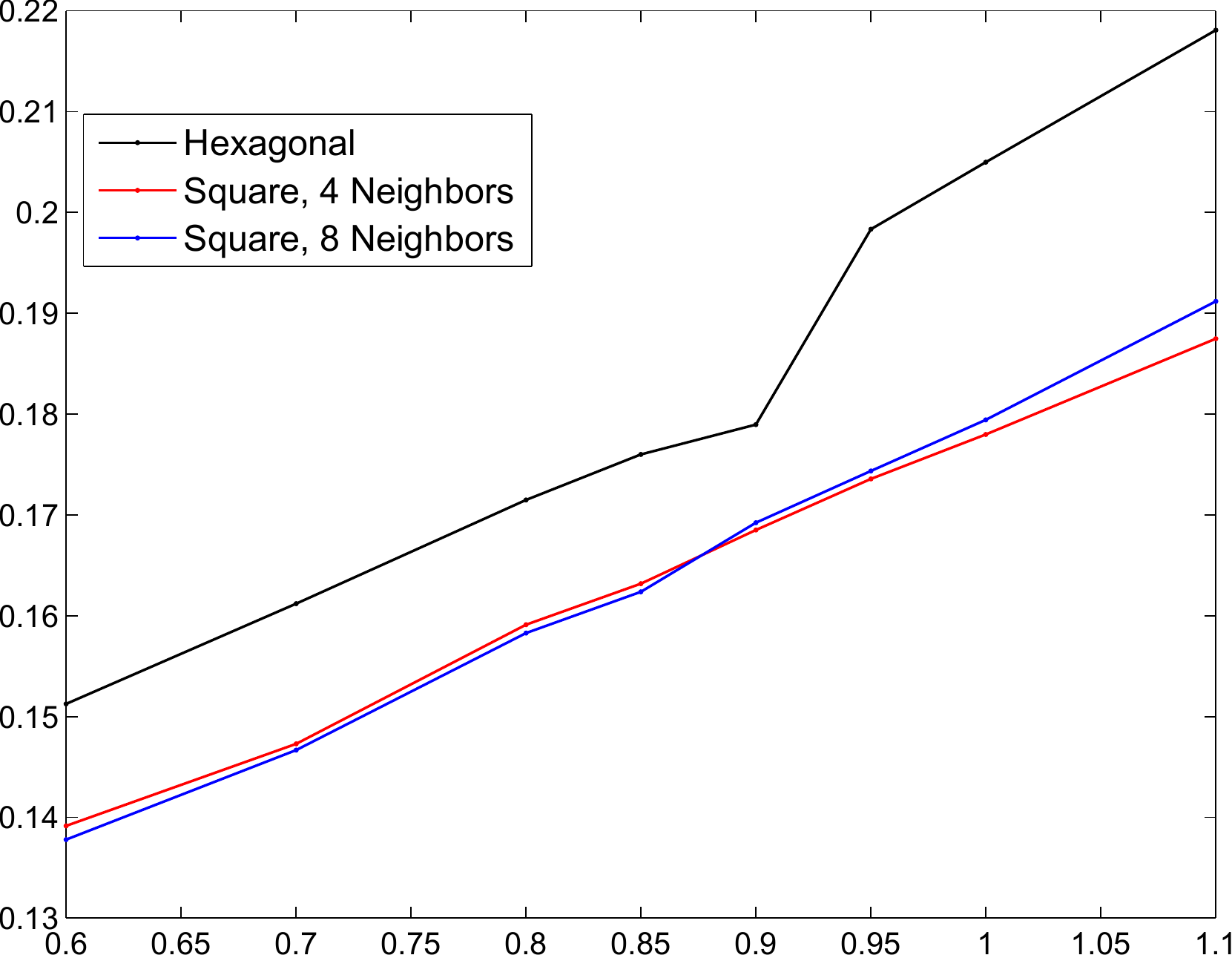}
  \caption{Results of Experiment No.~5. Averaged $\ell^1$ distances (left) and averaged ratios of correctly restored pixels (right) plotted against $\lambda$.}
  \label{fig:exp5}
\end{figure}

After adding $60\%$ salt \& pepper noise, the images were again denoised with an $L^1$ data term. Examples of denoised images are presented in Fig.~\ref{fig:exp4img}, error plots can be found in Fig.~\ref{fig:exp4}. The latter figure seems to indicate that the hexagonal structure achieved significantly better results for all reasonable values of $\lambda$. This has to be interpreted with care. It is possible, that the better performance figures of the hexagonal grid are caused by the grid conversion process simplifying the image.
\paragraph{Experiment No.\ 5}
Finally, the previous experiment was repeated with another natural image and $70\%$ salt \& pepper noise instead of $60\%$ (Figs.~\ref{fig:exp5img}, \ref{fig:exp5}).
\paragraph{Discussion}
A comparison of the error plots of all five experiments permits two remarkable conclusions. First, in every experiment TV denoising on the hexagonal grid achieved a smaller minimal $\ell^1$ distance between restored image and clean image. That is, the black curve reaches farthest down in Figs.\ \ref{fig:exp1}, \ref{fig:exp2}, \ref{fig:exp3}, \ref{fig:exp4} and \ref{fig:exp5}. In the majority of experiments the differences were, however, only marginal.

Second, if we denote by $\lambda^\star$ the distinguished value, where denoising on the hexagonal grid achieves the minimal $\ell^1$ error, then the hexagonal grid also performed better for all (depicted) $\lambda < \lambda^\star$. That is, the black curve lies below the others on the left hand side of Figs.\ \ref{fig:exp1}, \ref{fig:exp2}, \ref{fig:exp3}, \ref{fig:exp4} and \ref{fig:exp5}. This observation also holds true for the ratio of correctly restored pixels (with the only exception being Fig.~\ref{fig:exp2}). Note that those $\lambda$'s are more interesting, since larger values result in less filtered images, which at least from the viewpoint of visual appearance or restoration quality are rather futile (cf.\ Experiment No.\ 1 and also Fig.~\ref{fig:jump}).
\section{Conclusion}\label{sec:conclusion}
The problem of scalar image denoising by means of TV minimization under both $L^1$ and $L^2$ data fidelity was considered. We demonstrated how to reasonably transform this into a fully discrete (i.e.\ sampled and quantized) setting on an arbitrary lattice. Finally, by conducting a series of five experiments with both synthetic and natural images contaminated with different types of noise, we laid special emphasis on a comparison of TV regularization with respect to hexagonal and square discretizations.

With the caution in mind that, on the one hand, performance figures can never entirely capture the outcome of an experiment and that, on the other hand, evaluation by means of visual appearance is prone to subjectivity and misjudgement, it still seems safe to interpret our results as slightly favouring the hexagonal grid. This is not only due to the in general better restoration quality in terms of both $\ell^1$ errors and ratios of correctly restored pixels, but also due to its increased capability of suppressing unwanted metrication artefacts (and, seemingly, staircasing effects). Interestingly, the greatest discrepancies between the two discretization schemes occurred for the smoothest and most ``circular" of test images, the cosine in Experiment No.~2, which matches our initial motivation to study hexagonal grids.

Apart from giving comparative results, our work contributes to supplying the (admittedly still modest) demand for image processing tools handling non-standard images and, more importantly, makes a case for the significance of spatial discretization of continuous problems. Moreover, it suggests new research directions, such as generalization of the proposed discretization approach to higher dimensions of both image domain and range or an analysis of its approximation properties, which is planned to be tackled in future work.
\section*{Acknowledgements}
I take this opportunity to thank my colleagues at the Computational Science Center, University of Vienna, for the pleasure of working with them, and especially Otmar Scherzer for valuable discussions and thoughtful comments. This work has been supported by the Vienna Graduate School (\textit{Initiativkolleg}) in Computational Science funded by the University of Vienna as well as by the Austrian Science Fund (FWF) within the national research network Photoacoustic Imaging in Biology and Medicine, project S10505-N20.

\begin{singlespace}
\footnotesize
\def\cprime{$'$} \providecommand{\noopsort}[1]{}

\end{singlespace}

\end{document}